%% file: Triangle_decompositions_two_dimensional_simplicial_complexes.tex
\documentclass[11pt]{article}
\usepackage[margin=0.75in]{geometry}

\usepackage{amssymb,amsmath,mathtools,amsthm,thmtools,amsfonts}
\usepackage{thm-restate}

\usepackage{url} 
\usepackage{hyperref}
\hypersetup{
colorlinks,
linkcolor={red}, 
citecolor={black},
urlcolor={blue!60!black},
pdftitle={Triangle decompositions of two-dimensional simplicial complexes},
pdfauthor={Michelle Delcourt, Thomas Lesgourgues, Luke Postle}}

\usepackage[utf8]{inputenc}
\usepackage{color}
\usepackage{float} 
\usepackage{subcaption} 

\usepackage{enumitem}
\setlist[itemize]{topsep=.1in,itemsep=.1in}
\setlist[enumerate]{topsep=.1in,itemsep=.1in}
\setlist[enumerate,1]{label={(\roman*)}}
\setlist[enumerate,2]{label={(\alph*)}}
\setlist[enumerate,3]{label={(\arabic*)}}

\setlist[itemize]{nolistsep,noitemsep, topsep=0pt}

\newcommand{\bb}[1]{\ensuremath{\mathbb{#1}}}
\newcommand{\mc}[1]{\ensuremath{\mathcal{#1}}}%

\newcommand{\cH}{\mc{H}}

\newcommand{\cT}{\mc{T}}

\newcommand{\bR}{\bb{R}}

\DeclarePairedDelimiter{\parens}{(}{)}
\DeclarePairedDelimiter{\set}{\{}{\}}
\DeclarePairedDelimiter{\brackets}{[}{]}

\DeclarePairedDelimiter{\size}{\lvert}{\rvert}   

\usepackage[noabbrev,capitalise]{cleveref}

\theoremstyle{plain}
\newtheorem{thm}{Theorem}[section]
\newtheorem{lem}[thm]{Lemma}
\newtheorem{prop}[thm]{Proposition}

\newtheorem{cor}[thm]{Corollary}

\newtheorem{conj}[thm]{Conjecture}

\usepackage{etoolbox}
 \AtEndEnvironment{proof}{\setcounter{claim}{0}}

\theoremstyle{plain} 
\newcommand{\thistheoremname}{}
\newtheorem{genericthm}{\thistheoremname}

\theoremstyle{definition}
\newtheorem{definition}[thm]{Definition}

\crefname{equation}{}{}
\crefname{lem}{Lemma}{Lemmas}
\crefname{claim}{Claim}{Claims}
\crefname{thm}{Theorem}{Theorems}
\crefname{prop}{Proposition}{Propositions}
\crefname{enumi}{}{}
\crefname{lstlisting}{Code}{Codes}

\usepackage[dvipsnames, table]{xcolor}

\DeclareFontFamily{U}{mathx}{\hyphenchar\font45}
\DeclareFontShape{U}{mathx}{m}{n}{
      <5> <6> <7> <8> <9> <10>
      <10.95> <12> <14.4> <17.28> <20.74> <24.88>
      mathx10
      }{}
\DeclareSymbolFont{mathx}{U}{mathx}{m}{n}
\DeclareMathAccent{\widebar}{0}{mathx}{"73}

\title{Triangle decompositions of two-dimensional simplicial complexes}

\author{
Michelle Delcourt
\thanks{Department of Mathematics, Toronto Metropolitan University (formerly named Ryerson University),
Toronto, Ontario M5B 2K3, Canada {\tt mdelcourt@torontomu.ca}. Research supported by NSERC under Discovery Grant No. 2019-04269 and a Sloan Research Fellowship.} \and 
Thomas Lesgourgues
\thanks{School of Mathematics and Statistics, University of New South Wales, Australia. {\tt t.lesgourgues@unsw.edu.au}. Supported by Australian Research Council grant DP250101611.} \and
Luke Postle
\thanks{Combinatorics and Optimization Department,
University of Waterloo, Waterloo, Ontario N2L 3G1, Canada {\tt lpostle@uwaterloo.ca}. Partially supported by NSERC
under Discovery Grant No. 2026-04610.}}
\date{\today}

\input{tikz}

\begin{document}

\maketitle
\begin{abstract}
We study fractional triangle decompositions of two-dimensional simplicial complexes, unifying the settings of Nash-Williams' conjecture on triangle decompositions and the high-dimensional Dirac problem for Steiner triple systems.

We prove that every sufficiently large two-dimensional simplicial complex with minimum degree at least $(1-d_1)n$ and maximum missing codegree at most $d_2n$ admits a fractional triangle decomposition, provided $d_1/0.21+d_2/0.27\leq 1$. As an immediate consequence, we prove the asymptotic form of Lee's high-dimensional Dirac conjecture: every sufficiently large admissible $3$-uniform hypergraph with minimum codegree at least $(3/4+\varepsilon)n$ contains a Steiner triple system.

We also provide constructions establishing lower bounds on this fractional decomposition threshold problem.
\end{abstract}


\section{Introduction}

Graph decomposition problems are a central theme in extremal design theory, of which triangle cases are foundational instances. Given a graph, we are interested in knowing of we can partition its edge set into triangles. 
Such a partition, called a {\em triangle-decomposition}, works on two layers; the graph whose edges are to be partitioned, and the triangles that are allowed to be used in the decomposition. Let $G$ be a graph, denote by $\mc{T}(G)$ the set of all triangles in $G$, and let $\mc{H}\subseteq\mc{T}(G)$. An \mc{H}-decomposition of $G$ is a family of triangles from $\mc{H}$ that partitions $E(G)$. Its fractional analogue is a function $\omega:\mc{H}\to[0,1]$ such that, for every edge $e\in E(G)$,
\[\sum_{T\in\mc{H}\colon e\subseteq T}\omega(T)=1.\]
Thus the underlying object is a two-dimensional simplicial complex whose $1$-skeleton is $G$ and whose $2$-faces are the members of $\mc{H}$. Simplicial complexes have long been used to formulate problems in extremal combinatorics; Keevash and Mycroft~\cite{KM2015} studied matchings in hypergraph by asking in how many ways each face of a simplicial complex can be extended by one vertex; in their seminal work on combinatorial designs, Glock, K\"{u}hn, Lo and Osthus~\cite[Definition~4.6]{GKLO16} defined graph decompositions based on simplicial complexes, that is, decomposition of a graph into subgraphs from a prescribed family. We study the fractional version of this problem and ask how many triangles through each edge may be forbidden, that is, in $\mc{T}(G)\setminus\mc{H}$, given a lower bound on the minimum degree of the graph $G$. 
This problem unifies two very natural veins of research.  On one hand, the case when $\mc{H} = \mc{T}(G)$ corresponds is the subject of Nash-Williams' conjecture~\cite{nash1970unsolved} and has a long history. On the other hand, the case when $G$ is complete has garnered recent attention.  We discuss the history of these two separate streams in subsections that follow before returning to the simplicial setting for our main result.

\subsection{Triangle decompositions of graphs}
The first boundary case is obtained by setting $\mc{H}=\mc{T}(G)$. A graph $G$ is said to be {\em $K_3$-divisible} if $3\mid e(G)$ and $d_G(v)$, the degree of $v$ in $G$, is even for every vertex $v\in V(G)$. Nash-Williams~\cite{nash1970unsolved} conjectured in 1970 that every sufficiently large $K_3$-divisible graph $G$ on $n$ vertices with minimum degree at least $3n/4$ admits a triangle-decomposition. The factor $3/4$ is best possible, as demonstrated by Graham in an editorial note to Nash-Williams' article. A central role in the history of the problem has been played by so-called {\em fractional} relaxations. Together with the seminal work of Haxell and R\"{o}dl~\cite{HR01} (a regularity-based transformation from fractional to {\em approximate} decompositions) the iterative absorption breakthrough of Barber, K\"uhn, Lo, and Osthus~\cite{BKLO16}, shows that, for $c\geq3/4$, a fractional threshold $c$ yields exact
triangle-decompositions at every minimum degree $(c+\varepsilon)n$, for arbitrarily small $\varepsilon>0$. This made the determination of the fractional threshold the key ingredient to prove an asymptotic version of Nash-Williams' Conjecture. Many successive partial results improved the fractional threshold, from Garaschuk~\cite{garaschuk2014linear} in 2014, Dross~\cite{dross2016fractional} in 2016, Dukes and Horsley~\cite{dukes_minimum_2020} in 2020, and finally Delcourt and Postle~\cite{DP2021progress} who proved that a minimum-degree of at least $0.8274n$ suffices.

In a breakthrough result in 2026, Delcourt and Postle~\cite{DP26NW} resolved the conjecture in full. They first proved the fractional Nash-Williams' Conjecture; every large-enough graph $G$ on $n$ vertices with minimum degree at least $3n/4$ admits a fractional triangle-decomposition. Combined with the aforementioned results, this proves the asymptotic form of Nash-Williams' conjecture, that is, every large-enough $K_3$-divisible graph $G$ on $n$ vertices with minimum degree at least $(3/4+\varepsilon)n$ admits a triangle-decomposition.  

The guiding principle for their solution to the fractional problem is to reformulate the fractional decomposition problem as a linear program and use Farkas' Lemma to pass to its dual: if a fractional triangle-decomposition does not exist, then there is an edge-weighting which is non-negative on every triangle but has negative total weight. This idea already formed the basis of Garaschuk's thesis~\cite{garaschuk2014linear}. Delcourt and Postle attacked this dual obstruction problem with a remarkably effective discharging argument, transferring weight from positive to negative edges until the assumed negative total weight is contradicted. This combination of duality and discharging is a particularly useful concept, and form the main idea behind the present article. This approach has been used already in several works. Allsop, Bowtell, Lesgourgues, and Petrova~\cite{ABLP26} adapted the discharging argument to the partite setting arising from partial Latin square completion. Granet and Horsley~\cite{GH2026} determined the asymptotic $C_{\ell}$-decomposition threshold for every odd $\ell\geq 73$. Hou, Hu, and Wang~\cite{HHW2026}, and independently Granet, Horsley, and Santrio~\cite{GHS2026} used a quantitative discharging estimate originating in Delcourt and Postle's argument to settle the remaining odd lengths $5\leq \ell\leq 71$. Granet, Horsley, and Santrio also determined approximate decomposition thresholds for a range of other $3$-chromatic graphs. The dual discharging framework is therefore already emerging as a useful method beyond the original Nash-Williams problem. Delcourt and Postle then impressively removed the remaining $\varepsilon$ margin to prove an exact version of Nash-Williams' Conjecture. To do so, they introduced a stability approach, describing all near-extremal obstructions, and then combined this stability result with refined absorption.

\subsection{Steiner triple systems in dense 3-uniform hypergraphs}
The second aspect considered in the literature is the so-called high-dimensional Dirac's Theorem. We now fix $G:=K_n$, the complete graph on $n$ vertices, and let $\mc{H}\subsetneq\mc{T}(G)$. Observe that $K_n$ is $K_3$-divisible if and only if $n\equiv 1,3\mod{6}$. An $\mc{H}$-decomposition of $K_n$ is then a Steiner triple system contained in $\mc{H}$ (viewed as a $3$-uniform hypergraph), that is, a collection of triples covering every pair of vertices exactly once. As observed by Linial~\cite{linial2018challenges}, since a Steiner triple system partitions pairs rather than vertices, it is a natural higher-dimensional analogue of a perfect matching; this motivates the terminology \emph{high-dimensional Dirac problem}. The
relevant density parameter is the minimum codegree of $\mc{H}$, denoted here by $\delta_2(\mc{H})$, and defined by 
\[\delta_2(\mc{H}):=\min\set*{d_{\cH}(uv)\colon u,v\in \binom{V(\cH)}{2}},\]
where, for every pair of distinct vertices $u,v$ in $\mc{H}$, 
\[d_{\mc{H}}(uv)=\size*{\set*{w\in V(\cH)\colon uvw\in E(\cH)}}.\]
Lee~\cite{L25} conjectured the following.
\begin{conj}[Lee~\cite{L25}]\label{conj:MainConjSTS}
    There exists a constant $C$ such that the following holds for every sufficiently large integer $n \equiv 1,3 \mod{6}$.
    Let $\mc{H}$ be a $3$-uniform hypergraph on $n$-vertices with the minimum codegree at least $\frac{3}{4}n + C$. Then $\mc{H}$ contains a Steiner triple system.
\end{conj}
Lee constructed examples showing that the constant $3/4$ is asymptotically best possible, building a $3$-uniform hypergraph on $n$ with vertices with minimum codegree at least $\frac{3n}{4}-10$, and containing no Steiner triple system. Importantly, Lee proved an analogue of the reduction theorem of Barber, K\"{u}hn, Lo, and Osthus~\cite{BKLO16} in this setting. Formally, we define the shadow of a $3$-uniform hypergraph $\mc{H}$ to be the graph $\partial \mc{H}$ on $V(\mc{H})$ with edges the pairs contained in at least one hyperedge of $\mc{H}$, and its essential codegree to be $\delta_2^{\rm ess}(\mc{H}) :=\min\set{d_{\mc{H}}(uv)\colon~uv\in E(\partial\mc{H})}$. For $\varepsilon\in[0,1)$, Lee defined $\theta_{\rm STS}^*(\varepsilon)$ as the infimum of all $c\in[0,1]$ such that every sufficiently large $3$-uniform hypergraph $\mc{H}$ with $\delta(\partial\mc{H})\geq (1-\varepsilon)(n-1)$ and $\delta_2^{\rm ess}(\mc{H})\geq c(n-2)$ admits a fractional decomposition of its shadow using triples of $\mc{H}$. We finally let, 
\begin{equation}\label{eq:DefLeeLimit}
\theta_{\rm STS}^{*} := \lim_{\varepsilon\to 0^+} \theta_{\rm STS}^*(\varepsilon). 
\end{equation}

Lee proved that, for every $\varepsilon>0$ and every large enough integer $n\equiv 1,3\mod{6}$, if $\mc{H}$ is an $n$-vertex $3$-uniform hypergraph such that
\[ \delta_2(\mc{H}) \geq \parens*{ 
\max\set*{ \theta_{\rm STS}^{*}, \frac34 } +\varepsilon}n,\]
then $\mc{H}$ contains a Steiner triple system. This implies that, if $\theta_{\rm STS}$ denotes the asymptotic minimum-codegree threshold for the existence of a Steiner triple system, then (the lower bound is also included in~\cite{L25})
\[ \frac34\leq \theta_{\rm STS} \leq \max\set*{ \theta_{\rm STS}^{*}, \frac34 }. \]

Thus, proving $\theta_{\rm  STS}^*\leq 3/4$ suffices to prove an asymptotic version of~\cref{conj:MainConjSTS}. It follows from the work of Glock, K\"{u}hn, Lo, and Osthus~\cite{GKLO16} that $\theta_{\rm STS}^*<1$. After important progress by Lee~\cite{L25}, who showed $\theta_{\rm STS}^*(\varepsilon)\leq \frac{3 + \sqrt{57}}{12} < 0.88$ for all $\varepsilon>0$, the best known result to date was due to Zheng~\cite{Z24} with $\theta_{\rm STS}^*(\varepsilon) < 0.8579$ for all $\varepsilon>0$.\bigskip 

\subsection{The fractional-decomposition domain}
The Nash-Williams and high-dimensional Dirac problems are therefore two boundary cases of the same decomposition problem: in the former every triangle of $G$ is available, whereas in the latter
the underlying graph is complete. We study the general setting in which neither layer is complete, studying directly the two-dimensional simplicial complex problem. Let $(G,\mc{H})$ be a two-dimensional simplicial complex of order $n$ as defined in~\cref{def:2complex}, that is, $G\subseteq K_n$ and $\mc{H}\subseteq \cT(G)$. An integral triangle-decomposition selects $2$-faces (i.e. triangles in $\mc{H}$) covering each $1$-face (i.e. edges of $G$) exactly once, and a fractional triangle-decomposition assigns weights to $2$-faces so that each $1$-face receives total weight one. For every edge $e=uv\in E(G)$,  we define the {\em missing codegree of $uv$ in $\mc{H}$ with respect to $G$} as
\[d^{\rm miss}_{\mc{H},G}(uv):=\size*{N_G(u)\cap N_G(v)}-d_{\mc{H}}(uv).\]
We then obtain the maximum missing-codegree of $\mc{H}$ with respect to $G$ as 
\[\Delta^{\rm miss}_2(\mc{H},G):=\max\set*{d^{\rm miss}_{\mc{H},G}(uv)\colon uv\in E(G)}.\]

We are interested in the following. 

\begin{definition}[Fractional-decomposition domain]\label{def:OptimalDomain}
    Let $\mc{D}^*$ be the set of all pairs $(d_1,d_2)\in(0,1]\times[0,1]$ such that the following holds for every integer $n$ large enough. Let $(G,\mc{H})$ be a two-dimensional simplicial complex of order $n$ satisfying 
    \[\delta(G)\geq (1-d_1)n,\qquad \text{ and }\qquad \Delta^{\rm miss}_2(\mc{H},G)\leq d_2n.\]
    Then there exists a fractional triangle-decomposition of $(G,\mc{H})$, that is a fractional triangle-decomposition of  $G$ using triangles in $\mc{H}$.
\end{definition}
Although codegrees in an $n$-vertex $3$-uniform hypergraph are naturally normalized by $n-2$, throughout this paper we normalize the missing codegree by $n$. It avoids unnecessary second order issues  in the exact boundary constructions~\cref{sec:Lower}, and is in line with usual Nash-Williams type conjectures.

The restriction $d_1>0$ in~\cref{def:OptimalDomain} is deliberate. There exists no simple $n$-vertex graph satisfying $\delta(G)\geq n$. The high-dimensional boundary of $\mc{D}^*$ is therefore understood as the limit when $d_1\to0^+$, in line with Lee's definition of the limiting fractional Steiner triple-system threshold~\cref{eq:DefLeeLimit}.
The other boundary, when $d_2=0$, corresponds to the fractional Nash-Williams problem, as solved by Delcourt and Postle~\cite{DP26NW}, where every triangle of $G$ belongs to $\mc{H}$.

These two classical boundary cases yield immediate bounds on the possible domain $\mc{D}^*$. The work of Delcourt and Postle~\cite{DP26NW} on the fractional Nash-Williams threshold implies that $(d_1,0)\in\mc{D}^*$ for every $d_1\leq 1/4$, while the folklore lower bound excludes every point $(d_1,d_2)$ with $d_1>1/4$. Conversely, 
Zheng~\cite[Theorem~1.15]{Z24} proved that $\theta_{\rm STS}^*(\varepsilon)\leq 0.8579$ for every $\varepsilon\in[0,1)$. Let $x_*\approx 0.142165$ be the unique root in $[0,1/6]$ of $8x^3-22x^2+10x-1=0$. For every $2$-complex $(G,\mc{H})$ satisfying our $(d_1,d_2)$ assumptions and such that $2d_1+d_2\leq x_*$, for every edge $uv\in E(G)$, we have 
\[d_{\mc{H}}(uv) \geq |N_G(u)\cap N_G(v)|-d_2n \geq(1-2d_1-d_2)n\geq (1-x_*)n>0,\]
and therefore $uv\in \partial\mc{H}$. Conversely, since $\mc{H}\subseteq\mc{T}(G)$, every edge of $\partial\mc{H}$ is in $G$. Therefore $\partial\mc{H}=G$ and
\[ \delta^{\rm ess}_2(\mc{H}) = \min\set{d_{\mathcal H}(uv)\colon~uv\in E(G)} 
\geq (1-x_*)n.\]
Therefore Zheng’s result~\cite[Theorem~3.2]{Z24} yields a fractional triangle decomposition of $\partial\mc{H}=G$ using triangles in $\mc{H}$, hence,
\[ \set*{(d_1,d_2)\in(0,1]\times[0,1]\colon~2d_1+d_2\leq x_* } \subseteq\mc{D}^*.\]

In the other direction, Lee’s tripartite construction (see~\cref{prop:LowerLee,fig:LeeConst})
shows that every point $(d_1,d_2)$ with $d_2\geq1/3$  is not in $\mc{D}^*$. In particular the domain $\mc{D}^*$ from~\cref{def:OptimalDomain} must be contained in $(0,1/4]\times[0,1/3]$. A natural first step into determining $\mc{D}^*$ is the following.

\begin{definition}[Triangular Domain]\label{def:domain}
    For $(a,b)\in (0,1]\times [0,1]$, let 
    \[\mc{D}_{a,b}:=\set*{(d_1,d_2)\in(0,1]\times[0,1]\colon~d_1\in[0,a],~d_2\leq b-\frac{b}{a}d_1}.\]
    that is $\mc{D}_{a,b}$ is the subset of $[0,1]^2$ below the line joining the points $(0,b)$ and $(a,0)$.     
\end{definition}

We adapt the ``discharging in the dual'' framework from~\cite{DP26NW} to the two-dimensional simplicial complex setting and combine it with a simple first-order discharging rule to prove the following.

\begin{thm}\label{thm:MainSimplicial}
    Let $a=0.21$ and $b=0.27$. For every $(d_1,d_2)\in\mc{D}_{a,b}$, the following holds for every $n$ large enough. Let $(G,\mc{H})$ be a two-dimensional simplicial complex of order $n$ satisfying 
    \[\delta(G)\geq (1-d_1)n,\qquad \text{ and }\qquad \Delta^{\rm miss}_2(\mc{H},G)\leq d_2n.\]
    There exists a fractional triangle-decomposition of $(G,\mc{H})$, that is a fractional triangle-decomposition of  $G$ using triangles in $\mc{H}$.
\end{thm}


We did not try to optimize the constant $a,b$ from~\cref{thm:MainSimplicial}. A careful analysis of the discharging rule would allow $a$ to be increased slightly, to about $0.215$, while calculations with a ``non-uniform'' first order rule suggest values near $0.23$ are attainable. However, as noted by Delcourt and Postle in their resolution of Nash-Williams' Conjecture, such a first-order rule will not reach the optimal value of $a=1/4$; this would likely require a second-order (e.g. sponsorship) based argument. For our present article, the decisive feature is that we obtain $b>1/4$. 
\cref{thm:MainSimplicial} yields the following 
\begin{cor}
The fractional Steiner triple system threshold satisfies $\theta^*_{\rm STS}\le0.73.$    
\end{cor}
\begin{proof}
Let $a:=0.21$ and $b:=0.27$. Let $\varepsilon<t\in(0,a)$ be arbitrarily small reals, and define $s:=b\cdot(1-\frac{t}{a})$. Let $\mc{H}$ be an n-vertex 3-uniform hypergraph with shadow $G:=\partial\mc{H}$ satisfying
\[ \delta(G)\geq (1-\varepsilon)(n-1)
\qquad\text{ and }\qquad
\delta_2^{\rm ess}(\mc{H})\geq(1-s)(n-2).
\]
Since $t>\varepsilon$, we have $\delta(G)\geq(1-t)n$ for all sufficiently large integers $n$. Furthermore, every edge of $G$ has at most $n-2$ common neighbors, hence
\[ \Delta_2^{\rm miss}(\mc{H},G)
\leq (n-2)-\delta_2^{\rm ess}(\mc{H})
\leq s(n-2)
\leq sn.
\]
As \((t,s)\in\mathcal D_{a,b}\), by~\cref{thm:MainSimplicial} there exists a fractional triangle-decomposition of $G$ using triangles in $\mc{H}$. It follows that 
\[ \theta^*_{\rm STS}(\varepsilon)\leq 1-s =0.73+\frac97t,\]
and the result follows by taking $t$ (and $\varepsilon$) arbitrarily close to $0$.    
\end{proof}

Lee's reduction~\cite{L25} then implies that the asymptotic minimum codegree threshold for the existence of a Steiner triple system in a $3$-uniform hypergraph is exactly $3/4$, yielding the following asymptotic form of Conjecture~\ref{conj:MainConjSTS}.

\begin{cor}
    For every  positive real number $\varepsilon > 0$, the following holds for every sufficiently large integer $n\equiv 1,3 \mod{6}$.
    Every $n$-vertex $3$-uniform hypergraph with minimum codegree at least $\big(\frac{3}{4} + \varepsilon\big)n$ contains a Steiner triple system.
\end{cor}

\begin{figure}[H]
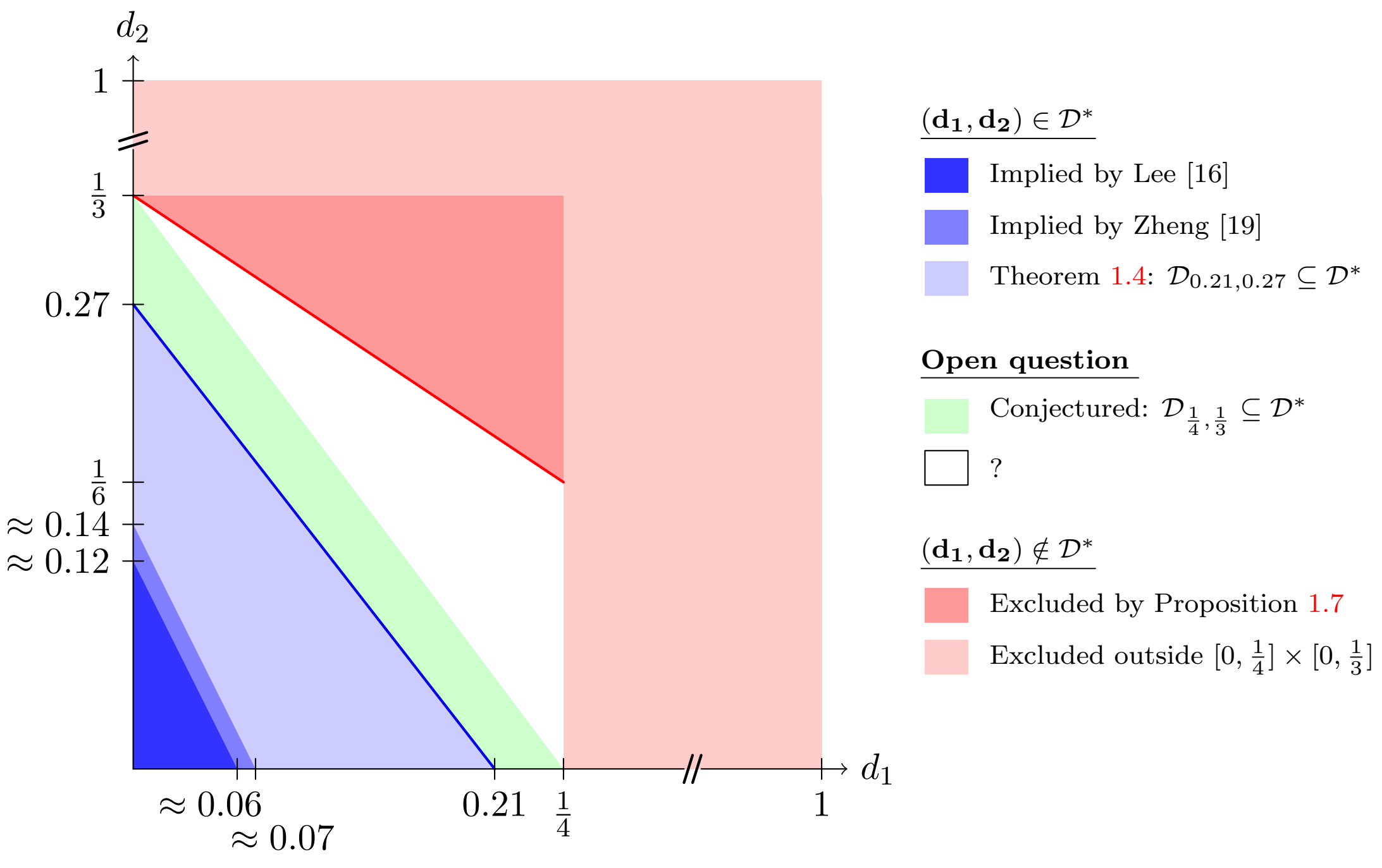

    \centering
    \DomainConjecture{18}{0.21}{0.27}
    \caption{The full picture: deciding if a given point $(d_1,d_2)$ is in $\mc{D}^*$.}
    \label{fig:conj}
\end{figure}

In addition to the positive result of~\cref{thm:MainSimplicial}, our second main result gives a non-trivial domain in $\parens*{\brackets*{0,\frac14}\times\brackets*{0,\frac13}}\setminus\mc{D}^*$. Our inner and outer bounds are summarized in
\cref{fig:conj}. We refer the reader to~\cref{sec:Conc} for a discussion on possible improvements and conjectures for the shape of $\mc{D}^*$.

\begin{prop}\label{lowerdomain}
    Let $\mc{Q}$ be defined as 
    \[\mc{Q}:=\set*{(d_1,d_2)\in[0,1]^2\colon d_1\leq\frac14,~2d_1+3d_2<1}.\]    We have 
    \[\mc{D}^*\subseteq \mc{Q}.\]
\end{prop}

\subsection{Proof overview}\label{sec:ProofOverview}

Our proof follows the ``discharging in the dual'' program of Delcourt and Postle~\cite{DP26NW}. Suppose that a $2$-complex $(G,\mc{H})$ admits no fractional triangle-decomposition. By Farkas' Lemma, there exists an edge-weighting
$\phi:E(G)\to\mathbb{R}$ of negative total value such that every triangle in $\mc{H}$ has a non-negative weight. 

We color an edge red if it has negative weight and blue otherwise. The proportional discharging rule from~\cref{def:Discharging} transfers weights from blue edges to adjacent red edges using only triangles in $\mc{H}$. In~\cref{prop:BlueStillPositive}, we show that blue edges never send more than their initial weight. We then fix a red edge $uv$ and in~\cref{prop:dischargingvalue}, we give a lower bound on the amount of weight it receives during the discharging. This produces an explicit function $f_{d_1,d_2}(x,y)$ considered in~\cref{prop:optim}.

The remaining task is analytic. We prove that $f_{d_1,d_2}(x,y)$ is at least $1$ for every admissible parameters $x,y$, as long as $(d_1,d_2)\in\mc{D}_{0.21,0.27}$. We obtain that every red edge receives at least the absolute value of its initial weight during discharging. Since the discharging principle ensures that the total weight is constant, this contradicts the negative total value of $\phi$.

\section{Preliminaries}

\begin{definition}\label{def:2complex}
A {\em two-dimensional simplicial complex}, or {\em $2$-complex for short}, is a pair $(G,\mc{H})$ where $G$ is a graph and $\mc{H}\subseteq\mc{T}(G)$. The vertices of $G$, edges of $G$, and members of $\mc{H}$ are respectively the $0$-, $1$-, and $2$-faces of the complex, and $G$ is its $1$-skeleton. The order of $(G,\mc{H})$ is $|V(G)|$.
\end{definition}

Following the work of Garaschuk~\cite{garaschuk2014linear} and of Delcourt and Postle~\cite{DP26NW}, we formulate the fractional $K_3$-decomposition problem as a linear program and use Farkas' Lemma to pass to the dual setting, where weights are assigned to the edges rather than the triangles. We first need some definitions, as introduced in~\cite{DP26NW}.

For a given $x\in \bR$, we let $x^+:=\max(x,0)$. We say that an edge-weight function (which we also refer to as an edge-weighting) $\phi:E(G)\to\mathbb{R}$ is {\em \mc{H}-covered} if every triangle $T\in\mc{H}\subseteq\cT(G)$ has non-negative total weight, that is,\[\sum_{e\in T}\phi(e)\geq 0,\text{ for all }T\in \cH.\]  
The {\em value} of $\phi$ is
\[|\phi|:=\sum_{e\in E(G)}\phi(e).\]
We say that $\phi$ is {\em totally-negative} if $|\phi|<0$.
	
We are now ready to state a version of Farkas' lemma adapted to our setting. This is the direct analogue of~\cite[Lemma~2.2]{DP26NW}, with the columns of edge-triangle incidence matrix being only the triangles contained in $\mc{H}$.

\begin{lem}[Farkas' Lemma applied to $K_3$-Decompositions]\label{lem:Farkas} A two-dimensional simplicial complex $(G,\mc{H})$ does not admit a fractional $K_3$-decomposition if and only if there exists an \mc{H}-covered and totally-negative edge-weighting of $G$. 
\end{lem}

Delcourt-Postle~\cite{DP26NW} also introduce the following concepts, which are useful property, make everything easier to work with.
We further say that an edge-weighting $\phi:E(G)\to\mathbb{R}$ is:
\begin{enumerate}[resume]
    \item {\em untied} if no two edges have the same weight, that is, 
    \(\phi(e)\neq \phi(f)\text{ for all }e\neq f\in E(G).\)
    \item {\em nonzero} if $\phi(e)\neq 0$ for all $e\in E(G)$.
\end{enumerate}

The following perturbation lemma is stated by Delcourt and Postle~\cite[Proposition~4.7]{DP26NW} when $\mc{H}=\mc{T}(G)$. Its proof uses only the inequalities $\sum_{e\in T}\phi(e)\geq 0$ for $T\in\mc{H}$,
and therefore applies verbatim to an arbitrary prescribed family
$\mc{H}\subseteq\mc{T}(G)$.

\begin{prop}[Delcourt-Postle~\cite{DP26NW}]\label{prop:NonZeroUntied}
Let $(G,\mc{H})$ be a two-dimensional simplicial complex.  If $\phi$ is a \mc{H}-covered edge-weighting of a graph $G$ with $e(G)\ge 1$ and $\varepsilon > 0$ is real, then there exists a nonzero untied \mc{H}-covered edge-weighting $\phi'$ of $G$ such that $\phi(e)\le \phi'(e)$ for all $e\in E(G)$ and 
\[\sum_{e\in E(G)} \phi(e)\leq \sum_{e\in E(G)} \phi'(e) \leq \varepsilon + \sum_{e\in E(G)} \phi(e).\]  
\end{prop}

We can then restate Farkas' Lemma with the following extension, whose proof is an immediate corollary of~\cref{prop:NonZeroUntied} by taking $\varepsilon<-|\phi|$.
\begin{lem}[Farkas’ Lemma extended]\label{lem:FarkasExtended} A two-dimensional simplicial complex $(G,\mc{H})$ does not admit a fractional $K_3$-decomposition if and only if there exists a nonzero, untied, \mc{H}-covered and totally-negative edge-weighting of $G$. 
\end{lem}

\section{First order discharging}

We now implement the proof strategy outlined in~\cref{sec:ProofOverview}. Fix a $2$-complex $(G,\mc{H})$ satisfying the assumptions of~\cref{thm:MainSimplicial}, and suppose for a contradiction that it admits no fractional triangle-decomposition. By~\cref{lem:FarkasExtended}, there exists a nonzero, untied, $\mc{H}$-covered and totally-negative edge-weighting $\phi$. We then define a discharging operation in~\cref{def:Discharging}, to transfer weights from positive (blue) edges to negative (red) edges, while preserving the total weight constant. We first prove in~\cref{prop:BlueStillPositive} that every blue edge has a non-negative final weight. 
For every red edge $r=uv$, we show in~\cref{prop:dischargingvalue} that it receives from blue edges at least an $f$-proportion of its initial weight, with $f:=f(d_1,d_2,x,y)$, for some $0<y\leq x<1$ with $x=x(u,v)$ and $y=y(u,v)$. We finally prove in~\cref{prop:optim} that, if $(d_1,d_2)\in\mc{D}_{a,b}$ from~\cref{def:domain}, then $f(d_1,d_2,x,y)\geq 1$ for all admissible values of $x,y$. This implies that each red edge gets at least its weights back during discharging, contradicting the fact that $\phi$ is totally-negative.


\subsection{Discharging rule}

Our argument uses the first-order discharging rule developed by Delcourt and Postle in their work on the fractional Nash-Williams conjecture. In the ordinary graph setting, they mention that these degree-based rule already show that minimum degree $0.785n$ suffices (while further first-order refinements reached $0.77n$). This simple rule is nevertheless particularly well suited to our $2$-complex setting; it yields a non-trivial inner bound for $\mc{D}^*$ and, at the high-dimensional boundary $d_1\to0^+$, it gets past the minimum codegree threshold of $3n/4$. Since this already suffices for our main application, while any first-order rule is unlikely to reach a complete determination of $\mc{D}^*$, we do not attempt to optimize the constants further.


Following~\cite{DP26NW}, there are two types of triangles in $\mc{H}$ containing negative edges: those with one negative edge and two positive edges, which we call {\em acute},  and those with two negative edges and one positive edge, which we call {\em obtuse}; all other triangles are either irrelevant for our problem (they do not contain any negative edges) or are not in $\cH$. This follows from the fact that we are only considering edge-weightings that are $\mc{H}$-covered, that is, every triangle of $\mc{H}$ contains at least one positive edge.

For a given $2$-complex $(G,\mc{H})$ and a weight assignment $\phi$ to the edges in $G$, we define the {\em proportional discharging} of $\phi$ to be the following transfer of weights between edges of $G$, depending on the types of triangles they belong to. This rule is the first-order rule developed by Delcourt and Postle~\cite[Section~2.1.3]{DP26NW} and used in the partite setting in~\cite[Definition~4.1]{ABLP26}. The only difference in our setting comes from the restriction to triangles belonging to $\mc{H}$.

\begin{definition}[Proportional discharging.]\label{def:Discharging}
    Let $(G,\mc{H})$ be a two-dimensional simplicial complex, and let $\phi:E(G)\to\mathbb{R}$. For every edge $e\in E(G)$, color $e$ red if $\phi(e)<0$ and blue otherwise. Let $b=vw$ be a blue edge. The {\em proportional discharging of $b$} is the following transfer of weights from $b$ to its red neighbors. For each acute or obtuse triangle $T\in\mathcal{H}$ containing $b$ and the red edge $r=uv$, we transfer the weight ${\rm dis}_{\phi}(b,r)$ from $b$ to $r$ as 
    \[{\rm dis}_{\phi}(b,r):=
    \begin{cases}
        \dfrac{\phi(b)}{2\cdot d_{Red}(v)} &\text{ if $T$ is acute and in $\mc{H}$},\\\\
        \dfrac{-\phi(r)}{2\cdot d_{Red}(v)}+\dfrac{-\phi(r)}{2\cdot d_{Red}(w)} &\text{ if $T$ is obtuse  and in $\mc{H}$},
    \end{cases}\]
    where $d_{Red}(x)$ is the number of red edges incident with the vertex $x$ in $G$. For simplicity, we extend the definition of ${\rm dis}_{\phi}$ to all pairs of edges by setting ${\rm dis}_{\phi}(e_1,e_2):=0$ if $\phi(e_1)<0$, if $\phi(e_2)\geq0$, or if $e_1,e_2$ are not in a common triangle of $\mc{H}$.  For a given red edge $r\in G$, we define $\Delta_{\phi}(r)$ to be the total weight received by the edge $r$ during discharging, that is, 
    \[\Delta_{\phi}(r):=\sum_{b\in G}{\rm dis}_{\phi}(b,r).\]
\end{definition}

Following Delcourt and Postle~\cite{DP26NW}, we call the weight of an edge after all transfers its \emph{final charge}. For every $e\in E(G)$, define 
\[ {\rm ch}^{\phi}_{F}(e) := \phi(e) 
-\sum_{r\in E(G)}{\rm dis}_{\phi}(e,r)
+\sum_{b\in E(G)}{\rm dis}_{\phi}(b,e).
\]
Therefore, if $b$ is a blue edge and $r$ is a red edge, we obtain
\[
{\rm ch}^{\phi}_{F}(b) = \phi(b)-\sum_{r\in E(G)}{\rm dis}_{\phi}(b,r),\qquad
{\rm ch}^{\phi}_{F}(r) = \phi(r)+\Delta_{\phi}(r).
\]
Importantly, as discharging only transfers weights between edges, we have
\[ \sum_{e\in E(G)}{\rm ch}^{\phi}_{F}(e)=|\phi|. \]

\begin{table}[ht]
    \centering
    \begin{tabular}{>{\centering\arraybackslash}p{0.4\textwidth} | >{\centering\arraybackslash}p{0.4\textwidth}}
    \textbf{$T$ is acute in $\mc{H}$} & \textbf{$T$ is obtuse in $\mc{H}$}\\\hline
    \DischargingQuasiPos{1.5} & \DischargingQuasiNeg{1.5} \\
    ${\rm dis}_{\phi}(b,r):=\dfrac{\phi(b)}{2\cdot d_{Red}(v)}$     & ${\rm dis}_{\phi}(b,r_i):=\dfrac{-\phi(r_i)}{2\cdot d_{Red}(v)}+\dfrac{-\phi(r_i)}{2\cdot d_{Red}(w)}$
    \end{tabular}
    \caption{Discharging rules from a positive edge $b$ to one of its negative neighbors.}
    \label{tab:Discharging}
\end{table}


We first prove that no blue edge sends more charge than its initially one, and hence every blue edge has non-negative final charge. 

\begin{prop}\label{prop:BlueStillPositive}
     Let $(G,\mc{H})$ be a two-dimensional simplicial complex, and let $\phi:E(G)\to\mathbb{R}$ be an \mc{H}-covered edge-weighting of $G$. For every edge $b$ such that $\phi(b)\geq 0$, we have
    \[{\rm ch}_F^{\phi}(b)\geq 0\]
\end{prop}

The proof of~\cref{prop:BlueStillPositive} is identical to~\cite[Proposition~4.2]{ABLP26}, with transfers limited to triangles in $\mc{H}$. We include the short argument for completeness.

\begin{proof}
Let $b=vw \in E(G)$ be such that $\phi(b)\geq 0$, and without loss of generality assume that $d_{Red}(v)\geq d_{Red}(w)$. If $d_{Red}(v)+d_{Red}(w)=0$, then $b$ is in no obtuse and no acute triangles, hence $\sum_{r\in E(G)}{\rm dis}^f_{\phi}(b,r)=0$. If $d_{Red}(v)>d_{Red}(w)=0$, then $b$ is contained only in acute triangles of $\mc{H}$ and all edges that \emph{\mc{H}-incident} to $b$ are of the form $uv$, hence
\[\sum_{r\in E(G)}{\rm dis}_{\phi}(b,r)\leq 
d_{Red}(v)\cdot\frac{\phi(b)}{2d_{Red}(v)}\leq \phi(b).\]
Assume now that $d_{Red}(v)\geq d_{Red}(w)>0$. Denote by $\mc{H}^-_b$ the family of obtuse triangles of $\mc{H}$ containing $b$. We have
\begin{align*}
    \sum_{r\in E(G)}{\rm dis}_{\phi}(b,r)
 \leq 
 &\sum_{u\in V(G)\colon uvw\in \mc{H}^-_b} \brackets*{
 \frac{|\phi(uv)|}{2d_v}+\frac{|\phi(uv)|}{2d_w}
 +\frac{|\phi(uw)|}{2d_v}+\frac{|\phi(uw)|}{2d_w}
 }\\
 &+ \parens*{d_v-|\mc{H}^-_b|}\cdot \frac{\phi(vw)}{2d_v} 
 + \parens*{d_w-|\mc{H}^-_b|}\cdot \frac{\phi(vw)}{2d_w}\\
 \leq &\sum_{u\in U\colon uvw\in \mc{H}^-_b} \brackets*{
 \frac{\phi(vw)}{2d_v}+\frac{\phi(vw)}{2d_w}
 }
 + \parens*{d_v-|\mc{H}^-_b|}\cdot \frac{\phi(vw)}{2d_v}
 + \parens*{d_w-|\mc{H}^-_b|}\cdot \frac{\phi(vw)}{2d_w} \\
 =&~\phi(vw),
\end{align*}
as desired, where for the second inequality we have used the fact that $\phi$ is triangle-covered, and where the first inequality is not an equality as there might exist some triangles containing $b$ not in $\mc{H}$.
\end{proof}

We now lower bound the total charge received by a given red edge depending on the red degrees of its vertices. This follows the same broad strategy as the partite argument used for partial Latin square completion in~\cite{ABLP26}, although the present non-partite setting and the restricted family $\mc{H}$ lead to a different optimization problem.

\begin{prop}\label{prop:dischargingvalue}
    Let $(d_1,d_2)\in (0,1]\times[0,1]$, and let $(G,\mc{H})$ be a two-dimensional simplicial complex of order $n$, such that $\delta(G)\geq (1-d_1)n$ and $\Delta^{\rm miss}_2(\cH,G)\leq d_2 n$. If there exists an assignment of weights $\phi:E(G)\to \mathbb{R}$ that is \mc{H}-covered and nonzero, then, for every edge $r=uv\in G$ with $\phi(r)<0$ and $d_{Red}(u)\geq d_{Red}(v)$, we have
    \[\frac{\Delta_{\phi}(r)}{|\phi(r)|}
    \geq \Big(\frac{1}{2(1+d_2-x)} + \frac{1}{2y}- \frac{1}{2x}\Big)\cdot(x-d_2-d_1)^+ 
        + \frac{(y-d_2-d_1)^+}{2(1+d_2-y)}+ \frac{1-d_2-2d_1}{2x}\]
    where $x:=d_{Red}(u)/n$, $y:=d_{Red}(v)/n$ satisfy $0<y\leq x<1$ and $0<x+y\leq 1+d_2$.
\end{prop}

\begin{proof}
For simplicity, we extend the weight assignment by setting $\phi(e)=0$ for every edge $e\in K_n\setminus G$. As every triangle $T\in\mc{H}$ contains no edge $e\in K_n\setminus G$, we still have $ \sum_{e\in T}\phi(e)\geq 0$ for every $T\in\mc{H}$. We now color every edge $e$ of $K_n$ blue if $\phi(e)\geq 0$, and red if $\phi(e)< 0$. let 
    \[ R = \{e\in E(G)\colon \phi(e)<0\},\text{ and }\quad B = \{e\in E(G)\colon \phi(e)\geq 0\}.\]
As $\sum_{e\in T}\phi(e)\geq 0$ for every $T\in\mc{H}$, any fully red triangle $T$ in $G$ is not in $\mc{H}$. In particular, any red edge $uv$ is in at most $d_2n$ triangles in $\cT(G)\setminus\mc{H}$. 

We now use the proportional discharging of~\cref{def:Discharging} to transfer weights between edges. Fix a red edge $r=uv\in R$, hence with $\phi(r)<0$, such that $d_{Red}(u)\geq d_{Red}(v)$. For simplicity, denote $d_u:=d_{Red}(u)$ and $d_v:=d_{Red}(v)$. Recall that $\Delta_{\phi}(r)$ is the sum of weights received by $r$ during the discharging, that is
\[\Delta_{\phi}(r):=\sum_{b\in B}{\rm dis}_{\phi}(b,r).\]
We first observe that $u,v$ have at most $d_2n$ common red-neighbors, hence $(d_u-1)+(d_v-1)\leq n-2 + d_2n$ hence
\begin{equation}\label{eq:conditionxy}
    d_u+d_v\leq (1+d_2)n.
\end{equation}

We partition the vertices of $G\setminus{\{u,v\}}$ into $5$ sets of vertices, depending on the red-neighborhood of $u$ and $v$, and on the family $\mc{H}$ (see~\cref{fig:RedNeighbors}):
\begin{enumerate}
    \item $R_u$ contains the common-neighbors $w$ of $u,v$ such that $uvw\in\mc{H}$ and $uw$ is red, that is,
    \[R_u:=\big\{w\colon uvw\in\mc{H},~uw\in R,~vw\in B\big\}.\]
    
    \item $R_v$ contains the common-neighbors $w$ of $u,v$ such that $uvw\in\mc{H}$ and $vw$ is red, that is,
    \[R_v:=\big\{w\colon uvw\in\mc{H},~vw\in R,~uw\in B\big\}.\]    
    
    \item $B_{uv}$ contains the common-neighbors $w$ of $u,v$ such that $uvw\in\mc{H}$, and $uw$ and $vw$ are both blue, that is,
    \[B_{uv}:=\big\{w\colon uvw\in\mc{H},~uw,vw\in B\big\}.\]

    \item $N_{G}$ contains vertices of $G$ that are a non-neighbor with $u$ and/or $v$, that is,
    \[N_{G}:=\big\{w\colon uvw\not\in \cT(G)\}.\]

    \item Finally, $N_{\mc{H}}$ contains the common-neighbors $w$ of $u,v$ such that $uvw\not\in\mc{H}$, that is, 
    \[N_{\mc{H}}:=\big\{w\colon uvw\in \cT(G)\setminus\mc{H}\}.\]
    
\end{enumerate}

\begin{figure}[ht]
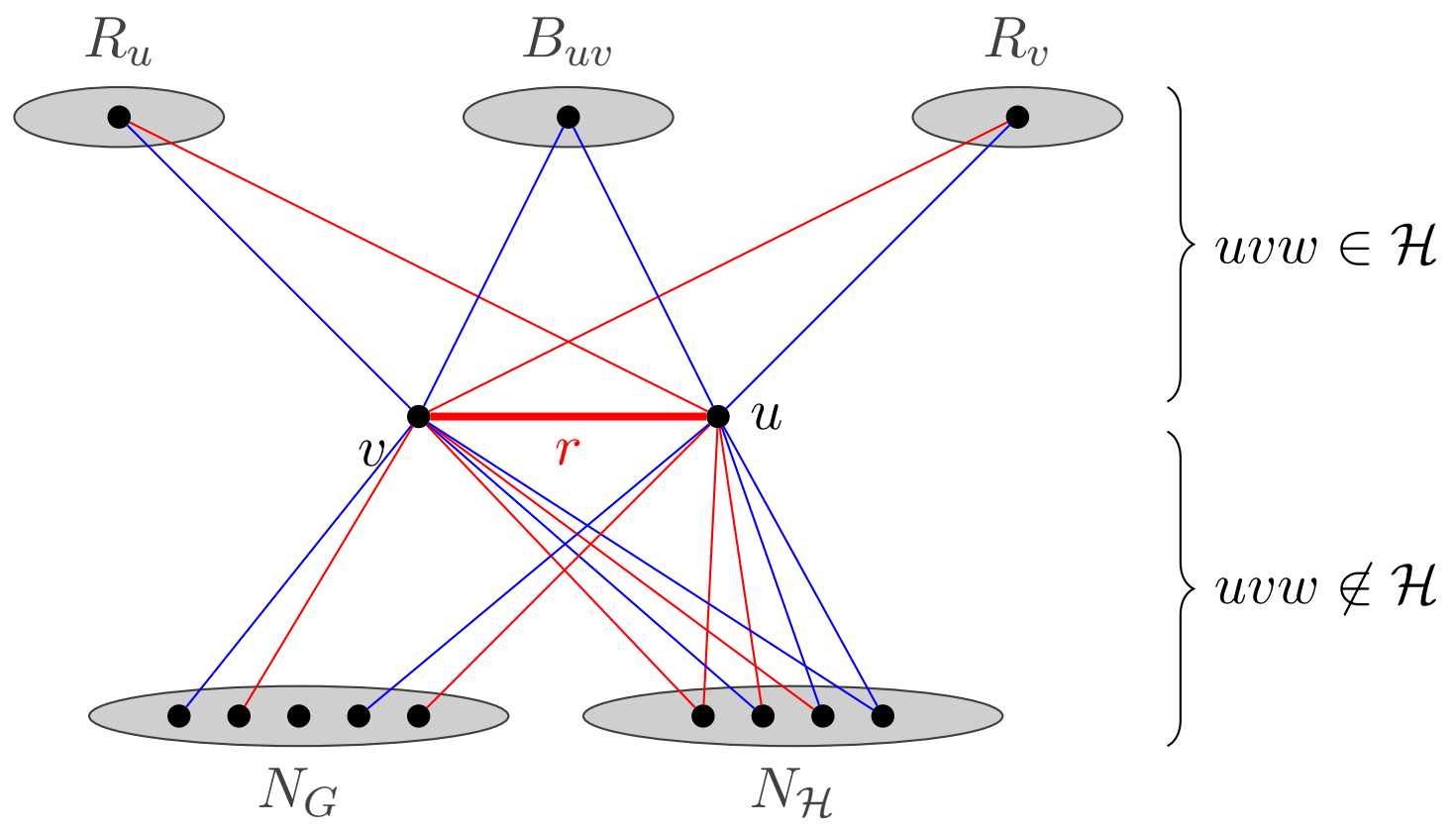

    \centering
    \RedEdges{2}
    \caption{Neighbors of $u,v$}
    \label{fig:RedNeighbors}
\end{figure}

We further define $N:=N_G\cup N_{\mc{H}}$. By definition of $(G,\cH)$, we have that 
\begin{equation}\label{eq:N}
|N_G|\leq 2(d_1n-1),\quad |N_{\cH}|\leq d_2n.    
\end{equation}

Observe that, for every $w \in V(G)$ such that $uw$ is red, then there are three possibilities: $vw$ is a blue edge, $vw\notin E(G)$, or $uvw\not\in\mc{H}$. Hence we have
\begin{align}
|R_u|   &\geq \big(d_u - 1 - (d_1n-1)-d_2n\big)^+,\label{eq:ObtuseU}\\
|R_v|   &\geq \big(d_v - 1 - (d_1n-1)-d_2n\big)^+,\label{eq:ObtuseV}\\
|B_{uv}| &= n-2-|R_u|-|R_v|-|N|. \label{eq:Acute}
\end{align}

\noindent For every vertex $w\in B_{uv}$, hence such that $T=uvw$ is an acute triangle in $\mc{H}$, through discharging the edge $r$ receives an additional weight of 
\[\frac{\phi(uw)}{2d_u}+\frac{\phi(vw)}{2d_v}\geq \frac{\phi(uw)+\phi(vw)}{2d_u}\geq \frac{|\phi(r)|}{2d_u},\]
where the last inequality follows from $\sum_{f\in T}\phi(f)\geq 0$. We now fix a vertex $w\in R_u$, hence such that $uvw$ is an obtuse triangle in $\mc{H}$ with $uw\in R$. Observe that any red-neighbor of $w$ that is also a red-neighbor of $u$ induces a negative triangle. As each such triangle cannot be in $\mc{H}$, we obtain that $w$ has at most $d_2n$ common red-neighbors with $u$, and therefore
\[d_{Red}(w)\leq (n-2)-(d_u-d_2n).\]

Moreover, as $uvw$ is an obtuse triangle with only blue edge $vw$, we have $\phi(vw)\geq |\phi(r)|$, and it follows that $r$ receives from $wv$ a weight of
\[ \frac{\phi(vw)}{2 d_w}+\frac{\phi(vw)}{2 d_v} \geq 
\frac{|\phi(r)|}{2\cdot (n+d_2n-d_u-2)}+\frac{|\phi(r)|}{2 d_v}.\]
We obtain that the edge $r$ receives, from the discharging process, a total weight of $\Delta_{\phi}(r)$ with
\begin{align*}
\frac{\Delta_{\phi}(r)}{|\phi(r)|}\geq 
    \parens*{\frac{1}{2\cdot (n+d_2n-d_u-2)}+\frac{1}{2d_v}}\cdot|R_u|
    +\parens*{\frac{1}{2\cdot (n+d_2n-d_v-2)}+\frac{1}{2d_u}}\cdot|R_v|
    +\frac{1}{2d_u}\cdot |B_{uv}|.
\end{align*}
Using~\cref{eq:Acute}, we obtain
\begin{align*}
\frac{\Delta_{\phi}(r)}{|\phi(r)|}\geq 
    \parens*{\frac{1}{2\cdot (n+d_2n-d_u-2)}+\frac{1}{2d_v}-\frac{1}{2d_u}}\cdot|R_u|
    +\frac{1}{2\cdot (n+d_2n-d_v-2)}\cdot|R_v|
    +\frac{1}{2d_u}\cdot (n-2-|N|).
\end{align*}

Observe that $\frac{1}{2d_v}-\frac{1}{2d_u}\geq 0$, therefore, using~\crefrange{eq:N}{eq:ObtuseV} we obtain
\begin{align*}
\frac{\Delta_{\phi}(r)}{|\phi(r)|}\geq 
     &\parens*{\frac{1}{2\cdot (n+d_2n-d_u)}+\frac{1}{2d_v}-\frac{1}{2d_u}}\cdot\big(d_u - (d_1+d_2)n\big)^+\\
    +&\frac{1}{2\cdot (n+d_2n-d_u)}\cdot\big(d_v - (d_1+d_2)n\big)^+\\
    +&\frac{1}{2d_u}\cdot \Big(n-(2d_1+d_2)n\Big).
\end{align*}

We normalize by $n$ by defining $x:=\frac{d_u}{n}$, and $y:=\frac{d_v}{n}$, with $0<y\leq x< 1$. We obtain
\begin{align*}
    \frac{\Delta_{\phi}(r)}{|\phi(r)|}&\geq 
        \parens*{\frac{1}{2\cdot (1+d_2-x)}+\frac{1}{2y}-\frac{1}{2x}}\cdot(x-d_1-d_2)^+
        +\frac{1}{2\cdot (1+d_2-y)}\cdot(y-d_1-d_2)^+
        +\frac{1-2d_1-d_2}{2x}.
\end{align*}
It also follows from~\cref{eq:conditionxy} that $x+y\leq 1+d_2$, as desired.
\end{proof}

\subsection{Optimization}

\begin{prop}\label{prop:optim}
    Let $a=0.21$ and $b=0.27$. Let $f_{d_1,d_2}:[0,1]^2\to\bR$ be defined as 
    \[f_{d_1,d_2}(x,y)=\Big(\frac{1}{2(1+d_2-x)} + \frac{1}{2y}- \frac{1}{2x}\Big)\cdot(x-d_2-d_1)^+ 
        + \frac{(y-d_2-d_1)^+}{2(1+d_2-y)}+ \frac{1-d_2-2d_1}{2x}.\]
    Then for every $(d_1,d_2)\in\mc{D}_{a,b}$ as in~\cref{def:domain}, we have $f_{d_1,d_2}(x,y)\geq 1$ for all $x,y\in(0,1)$ with $x\geq y$.
\end{prop}

\begin{proof}
By definition of $\mc{D}_{a,b}$, we have $d_2\leq b-\frac{b}{a}d_1$ hence $d_2\leq b-\frac97d_1$. Using the $b=0.27<\frac27<\frac13$ we obtain that 
\begin{equation}\label{eq:sum}
    d_1+d_2\leq b-\frac{2}{7}d_1\leq b(1-d_1)<\frac13(1-d_1).
\end{equation}

Assume first that $y\leq x \leq d_1+d_2$. Using~\cref{eq:sum} we obtain
\[f_{d_1,d_2}(x,y) =  \frac{1-d_2-2d_1}{2x}\geq \frac{2d_1+2d_2}{2x}\geq 1,\]
as desired. Assume now that $d_1+d_2\leq x$. For simplicity, let $d:=d_1+d_2$ hence such that $d\leq b(1-d_1)$ and $d\leq x$. By definition of $f_{d_1,d_2}$ we have
\[2f_{d_1,d_2}(x,y) = \frac{x-d}{1+d_2-x}+\frac{x-d}{y}+\frac{(y-d)^+}{1+d_2-y}+\frac{1-x-d_1}{x}.\]
Observe that
\[1+\frac{(y-d)^+}{1+d_2-y}=\frac{1-d_1}{1-d_1-(y-d)^+},\]
therefore we obtain
\[2f_{d_1,d_2}(x,y) +2= 
\frac{x-d}{1+d_2-x}+\frac{x-d}{y}+\frac{1-d_1}{1-d_1-(y-d)^+}+\frac{1-d_1}{x}.\]
By the Engel form of Cauchy–Schwarz inequality ($\frac{a^2}{u}+\frac{b^2}{v}\ge \frac{(a+b)^2}{u+v}$ for $u,v>0$), we have
\[\frac{x-d}{y}+\frac{1-d_1}{1-d_1-(y-d)^+}\geq
\frac{(\sqrt{x-d}+\sqrt{1-d_1})^2}{y+1-d_1-(y-d)^+}
\geq \frac{(\sqrt{x-d}+\sqrt{1-d_1})^2}{1-d_1+d}
= \frac{(\sqrt{x-d}+\sqrt{1-d_1})^2}{1+d_2},\]
hence
\[2f_{d_1,d_2}(x,y) +2\geq 
\frac{x-d}{1+d_2-x}+\frac{1-d_1}{x} + \frac{(\sqrt{x-d}+\sqrt{1-d_1})^2}{1+d_2}.\]
Let $r:=\sqrt{\frac{x-d}{1-d_1}}\in[0,1)$, hence such that
\[(1+d_2-x)=(1-d_1)(1-r^2),\quad x=(1-d_1)(r^2+\frac{d}{1-d_1}),
\quad (1+d_2)=(1-d_1)(1+\frac{d}{1-d_1}).\]
Using that $d\leq b(1-d_1)$ we obtain
\[(1+d_2-x)=(1-d_1)(1-r^2),\quad x\leq (1-d_1)(r^2+b),
\quad (1+d_2)\leq (1-d_1)(1+b),\]
and therefore
\[2f_{d_1,d_2}(x,y) +2\geq 
\frac{r^2}{1-r^2}+\frac{1}{r^2+b} + \frac{(r+1)^2}{1+b}:=g(r).\]

We have, for every $r\in[0,1)$,
\begin{align*}
g'(r)
&= \frac{2r(1-r^2)+r^2(2r)}{(1-r^2)^2}-\frac{2r}{(r^2+b)^2}+\frac{2(1+r)}{1+b}\\
&= \frac{2r}{(1-r^2)^2}-\frac{2r}{(r^2+b)^2}+\frac{2(1+r)}{1+b},\\
g''(r)
&=\frac{2(1-r^2)^2-2r(-4r(1-r^2))}{(1-r^2)^4} 
-\frac{2(r^2+b)^2-2r(4r)(r^2+b)}{(r^2+b)^4} +\frac{2}{1+b}.    \\
&=\frac{2(1+3r^2)}{(1-r^2)^3} +\frac{2(3r^2-b)}{(r^2+b)^3} +\frac{2}{1+b},\\    
\intertext{and,}
g'''(r) &= \frac{24r(1+r^2)}{(1-r^2)^4} +\frac{24r(b-r^2)}{(r^2+b)^4}.
\end{align*}

One can then check that, on the interval $[0,3/10]$, $g'''(r)>0$, hence $g''$ is increasing; We then have $g''(0)<0$ while $g''(3/10)>0$, hence $g'$ decreases and then increases; $g'(0)>0$ while $g'(3/10)<0$, hence $g'$ has exactly one zero on $[0,3/10]$, where $g$ is maximized. Therefore the minimum of $g$ on this interval is attained at an endpoint. We check that $g(0)>4$ and $g(3/10)>4$, therefore $g(r)>4$ for all $r\in[0,3/10]$.

Then, on the interval $[3/10,1]$, we have $g''(r)>0$, while $g'(3/10)<0$ and $g'(r)\to\infty$ as $r\to1$. Therefore $g$ has a unique minimum on $[3/10,1)$, and $g'$ is strictly increasing. We then have $g'(9/20)<0$ and $0<g'(1/2)<\frac49$, hence the minimizer $r_0\in[3/10,1]$ such that $g'(r_0)=0$ is in $(9/20,1/2)$. Then
    \[ g\parens*{\frac12}-g(r_0) =\int_{r_0}^{1/2}g'(t)dt 
    \leq \parens*{\frac12-r_0}\cdot g'\parens*{\frac12} 
    \leq \parens*{\frac12-\frac{9}{20}}\cdot g'\parens*{\frac12} <
    \frac{1}{20}\cdot\frac49
    =\frac1{45}.\]
As we also have $g(\frac12)=4+\frac{139}{4953}>4+\frac1{36}$, we obtain 
    \[g(r_0)>4+\frac{1}{36}-\frac{1}{45}>4,\]
    as desired.
\end{proof}

\subsection{Proof of Main Result}

We now complete the proof of~\cref{thm:MainSimplicial}. The discharging rule preserves the total weight of all edges;
\Cref{prop:BlueStillPositive} yields non-negative final charge for blue edges, while~\cref{prop:dischargingvalue,prop:optim} yield a similar result for red edges.

\begin{proof}
Fix $a=0.21$ and $b=0.27$, and let $(d_1,d_2)\in \mathcal{D}_{a,b}$. 
For $n$ large enough, let $(G,\mc{H})$ be a two-dimensional simplicial complex of order $n$, such that $\delta(G)\geq (1-d_1)n$ and $\Delta^{\rm miss}_2(\cH,G)\leq d_2 n$. 
Towards a contradiction, assume that $(G,\mc{H})$ admits no fractional triangle-decomposition. By~\cref{lem:FarkasExtended}, there exists a nonzero, untied, $\mc{H}$-covered and totally-negative edge-weighting $\phi$ of $G$.

Let $\psi:E(G)\to\mathbb{R}$ be the edge-weighting obtained from $\phi$ after applying the first-order discharging ${\rm dis}_{\phi}$ as defined in~\cref{def:Discharging}. By~\cref{prop:BlueStillPositive}, for every edge $b$ of $G$ with $\phi(b)\geq0$ we have
\[\psi(b) = \phi(b)-\sum_{r\in E(G)}{\rm dis}_{\phi}(b,r)\geq 0.\]
By~\cref{prop:dischargingvalue,prop:optim} for every edge $r\in G$ with $\phi(r)<0$, we have
\[\psi(r) = \phi(r)+\Delta_{\phi,r}=|\phi(r)|\cdot\parens*{\frac{\Delta_{\phi,r}}{|\phi(r)|}-1}\geq 0.\]
It follows that 
\[|\phi|=\sum_{e\in E(G)}\phi(e)=\sum_{e\in E(G)}\psi(e)\geq 0,\]
contradicting the fact that $\phi$ is totally-negative. Therefore $G$ admits a fractional $K_3$-decomposition.
\end{proof}

\section{Lower bounds}\label{sec:Lower}

The proof of~\cref{thm:MainSimplicial} yields an inner bound for $\mc{D}^*$. We now turn our attention to outer bounds. As mentioned before~\cref{def:domain}, the classical extremal constructions for the fractional Nash-Williams and fractional STS problems exclude the regions 
$d_1>1/4$ and $d_2\geq1/3$, respectively. However, these restrictions take no interaction into account. In particular, there is no reason to
expect points arbitrarily close to $(1/4,1/3)$ to belong to $\mc{D}^*$.

\subsection{Intuition}\label{sec:IntuitionLower}

We strongly believe that a more advanced discharging proof would improve the inner bounds on the domain $\mc{D}^*$, potentially up to 
\[\mc{D}_{\frac14,\frac13}=\set*{(d_1,d_2)\in(0,1]\times[0,1]\colon 4d_1+3d_2\leq1}.\]
However we do not believe that one to be the full picture. For example, we have some evidence that the restriction $\mc{D}^*\subseteq \mc{Q}$ should hold, where
\[\mc{Q}:=\set*{(d_1,d_2)\in[0,1]^2\colon d_1\leq\frac14,~2d_1+3d_2<1}.\]

Indeed, observe that, to the best of our knowledge, all known constructions for Nash-Williams type problems or higher-dimensional ones rely on building a graph $G$ with two types of edges, assigning a positive weight to $B\in E(G)$ and a negative one to $R\subseteq E(G)$, leading to a contradiction via Farkas Lemma. Therefore, assume that $(G,\mc{H})$ is a $2$-complex, with $E(G)=B\sqcup R$, $\delta(G)=(1-d_1)n$ and $\Delta^{\rm miss}_2(\mc{H},G)= d_2n$ for some $(d_1,d_2)\in(0,\frac14]\times[0,\frac13]$. Assume further that there exists an edge-weighting $\phi$ of $G$ such that $|\phi|<0$ and $\phi$ is $\mc{H}$-covered, and such that 
\[\phi(e)=\begin{cases}
+2&\text{ if }e\in B,\\
-1&\text{ if }e\in R.
\end{cases}\]
Therefore we have 
\[2 |B|-|R|<0,\qquad |R|>\frac23\cdot e(G)\geq\frac{1-d_1}{3}\cdot n^2.\]
As $d_1\leq 1/4$, we have $|R|\geq n^2/4$. Therefore, we know (see e.g. Bollob{\'a}s and Nikiforov~\cite[Corollary~2]{BN2005}, credited to Edwards~\cite{Edwards1977}, unpublished manuscript) that there exists an edge of $R$ that is contained in at least
\[\frac{2|R|}{n}-\frac{n}{3}\]
triangles fully in $R$. As $\phi$ is $\mc{H}$-covered, each such triangle is not in $\mc{H}$, and we obtain that 
\[d_2 \geq \frac{2|R|}{n^2}-\frac{1}{3}\geq \frac13-\frac{2d_1}{3},\]
that is 
\[\boxed{3d_2+2d_1\geq 1.}\]

In the opposite direction, for some integer $s,t$, let $G$ be a Tur\'{a}n graph with $s$ parts, each of size $3t$, that is the complete $s$-partite graph with parts $V_1,\ldots,V_s$ and $|V_1|=\ldots=V_s=3t$. Let $n:=3ts$. Further partition each part $V_i$ into three subparts $V_{i,1},V_{i,2},V_{i,3}$ each of size $t$. For every edge of $G$ let 
\[\phi(e)=\begin{cases}
+2&\text{ if }e\in E(V_{i,k},V_{j,k})\text{ for some } i,j\in[s], k\in[3],\\
-1&\text{ otherwise}.
\end{cases}\]
Observe that $|\phi|=0$. Delete one single positive edge from $G$, so that $|\phi|<0$. Let $\mc{H}$ be all triangle of $G$ that are not made of three negative edges. Any negative edges is in at most $(s-2)\cdot t$ negative triangles, therefore $\Delta^{\rm miss}_2(\mc{H},G)=d_2n$ and $\delta(G)=(1-d_1)n$ with 
\[d_2=\frac{s-2}{3s},\qquad d_1=\frac1s+\frac{1}{3ts}.\]
We obtain that
\[\boxed{3d_2+2d_1=1+\frac{2}{3ts}.}\]
Letting $s$ to be arbitrarily large, we obtain that the line $3d_2+2d_1=1$ is asymptotically sharp when $(d_1,d_2)$ tends to $(0,1/3)$.\bigskip

The equation $2d_1+3d_2=1$ suggests the uses of a 5-partite construction, 2 parts of size $d_1n$ and three of size $d_2n$. Our goal is to find one such construction $G$, with a family $\mc{H}\subseteq \mc{T}(G)$, both with small enough missing (co)degrees, with a totally-negative and $\mc{H}$-covered edge-weighting. Indeed, we have the following first example. 

\begin{prop}\label{prop:FirstPointNotD}
    The point $(d_1,d_2)=(\frac14,\frac16)$ is not in $\mc{D}^*$.
\end{prop}
\begin{proof}
    Let $n\equiv 0\mod{12}$. Let $G$ be the complete $5$-partite graph with parts $(U_1,U_2,U_3,V_1,V_2)$ such that $|U_1|=|U_2|=|U_3|=\frac{n}{6}$ and $|V_1|=|V_2|=\frac{n}{4}$, see~\cref{fig:FirstConst}. Let $\phi$ be the edge weighting of $G$ defined by
    \[\phi(uv)=
    \begin{cases}
        +1&\text{ if }u\in U_i,~v\in V_j\text{ for some }i\in[3],~j\in[2],\\
        -2&\text{ otherwise}.
    \end{cases}\]
    Observe first that
    \[|\phi|=(-2)\cdot \brackets*{3\cdot\parens*{\frac{n}{6}}^2+\parens*{\frac{n}{4}}^2}+ \parens*{\frac{n}{2}}^2=\frac{-1}{24}\cdot n^2<0,\]
    hence $\phi$ is totally-negative. We note then that the only negative triangles in $\mc{T}(G)$ are of the form $u_1u_2u_3$ with $u_i\in U_i$ for every $i\in [3]$ (all other triangles have null total weight). Therefore let 
    \[\mc{H}=\mc{T}(G)\setminus\set*{u_1u_2u_3\colon u_i\in U_i\text{ for every }i\in [3]}.\]
    We obtain that $\phi$ is $\mc{H}$-covered. Therefore by~\cref{lem:Farkas}, $(G,\mc{H})$ admits no fractional triangle-decomposition, with
    \[\delta(G)=\frac{3n}{4},\quad \Delta^{\rm miss}_2(\mc{H},G)=\frac{n}{6},\]
    as desired.
\end{proof}

\begin{figure}[H]
    \centering
    \FirstConstruction    
    \caption{Proving that $(\frac14,\frac16)\not\in\mc{D}^*$}
    \label{fig:FirstConst}
\end{figure}

\subsection{The 3-uniform lower bounds}
Before trying to generalize the construction from~\cref{prop:FirstPointNotD}, we now mention the construction used by Lee (and explained by Zheng in~\cite{Z24}) to prove that $\theta_{\rm STS}^* \geq 2/3$, which motivated his conjecture that $\theta_{\rm STS}^* = 2/3$.
\begin{prop}[Lee (see e.g.~\cite{Z24})]\label{prop:LowerLee}
	Consider the $3$-uniform hypergraph $\mc{H} = (V, E)$ on $n \geq 5$ vertices which is constructed as follows: Take an equitable partition $V_1 \cup V_2 \cup V_3$ of $V$ and set
	\[ E = V^{(3)} \setminus \set{e \in V^{(3)} \colon |e \cap V_1| = |e \cap V_2| = |e \cap V_3| = 1},\]
    that is $\mc{H}$ contains all triangle except for the transversal ones. Then $\delta_2(\mc{H})\geq 2n /3 - 8/3$ and $\mc{H}$ doesn't contain a fractional Steiner triple system.
\end{prop}

This is now an immediate consequence of~\cref{lem:Farkas}, by assigning a weight of $+2$ to every edge inside a part, and a weight of $-1$ to every edge between parts. Each triangle has total weight $0$ or $6$, while  the sum of weights over all edges is 
\[2\cdot 3\cdot\binom{n/3}{2} - 3\cdot \frac{n^2}{9}=-n<0.\]

\begin{figure}[H]
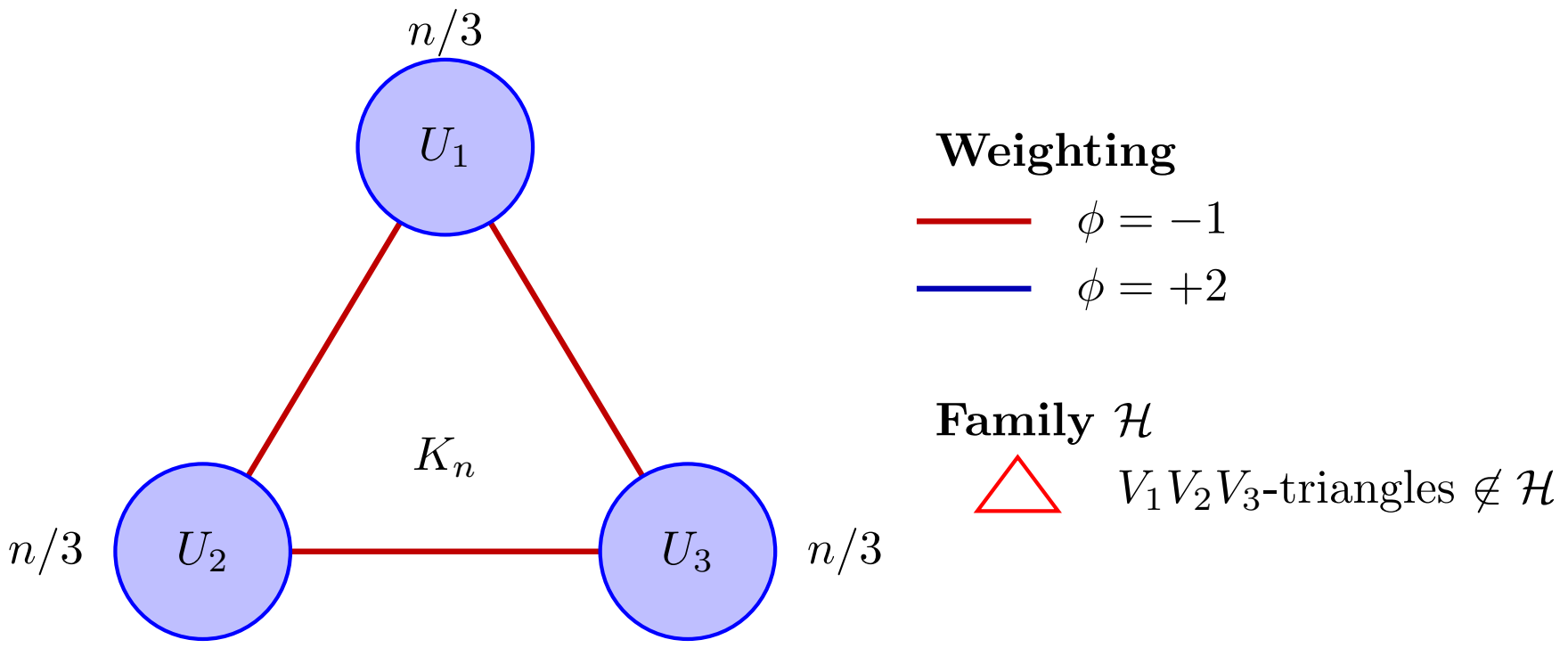

    \centering
    \LeeBis
    \caption{Lee's construction.}
    \label{fig:LeeConst}
\end{figure}

\subsection{Bounding the domain \texorpdfstring{$\mc{D}^*$}{D*}}

We now prove additional restriction on the domain $\mc{D}^*$,~\cref{lowerdomain}, based on the above two constructions,~\cref{fig:FirstConst,fig:LeeConst}. The main idea is a linear transition between the $3$-parts construction and the $5$-parts one. 
Recall that we defined \[\mc{Q}:=\set*{(d_1,d_2)\in[0,1]^2\colon d_1\leq\frac14,~2d_1+3d_2<1},\]
see~\cref{fig:conj}, and we now prove that $\mc{D}^*\subseteq \mc{Q}$.

\begin{proof}
Fix $d_1\in(0,\frac14]$, and let $d_2=\frac{1-2d_1}{3}$. As suggested by the equation $2d_1+3d_2=1$, we start with a complete $5$-partite graph on $n$ vertices, with 2 parts $V_1,V_2$ each of size $d_1n$, and three parts $U_1,U_2,U_3$ each of size $d_2n$. We want our new construction to match~\cref{fig:FirstConst} when $d_1=\frac14$, and to get arbitrarily close to~\cref{fig:LeeConst} when $d_1$ tends to $0$. In particular, we need to add edges inside the parts $U_i$'s ad $d_1$ decreases. Specifically, inside each part $U_i$, $i\in[3]$, we add a $(d_2-d_1)^+n$-regular graph. We obtain a graph $G$ with minimum degree exactly $(1-d_1)n$, as desired.

We now need to define an edge-weighting $\phi$ that matches both constructions. To do so, let $\phi:E(G)\to\mathbb{R}$ be defined by
\[\phi(uv)=
\begin{cases}
    +1 &\text{ if }u\in U_i,~v\in V_j\text{ for some }i\in[3],~j\in[2],\\
    +4 &\text{ if }u,v\in U_i\text{ for some }i\in[3],\\
    -2&\text{ otherwise}.
\end{cases}\]

As before, let $\mc{H}$ consist of all triangles of $G$ except
the negative ones, namely, except triangles $u_1u_2u_3$ with $u_i\in U_i$ for every $i\in [3]$. We obtain that $\Delta^{\rm miss}_2(\mc{H},G)=d_2n$, see~\cref{fig:InterpolConst}.
\begin{figure}[ht]
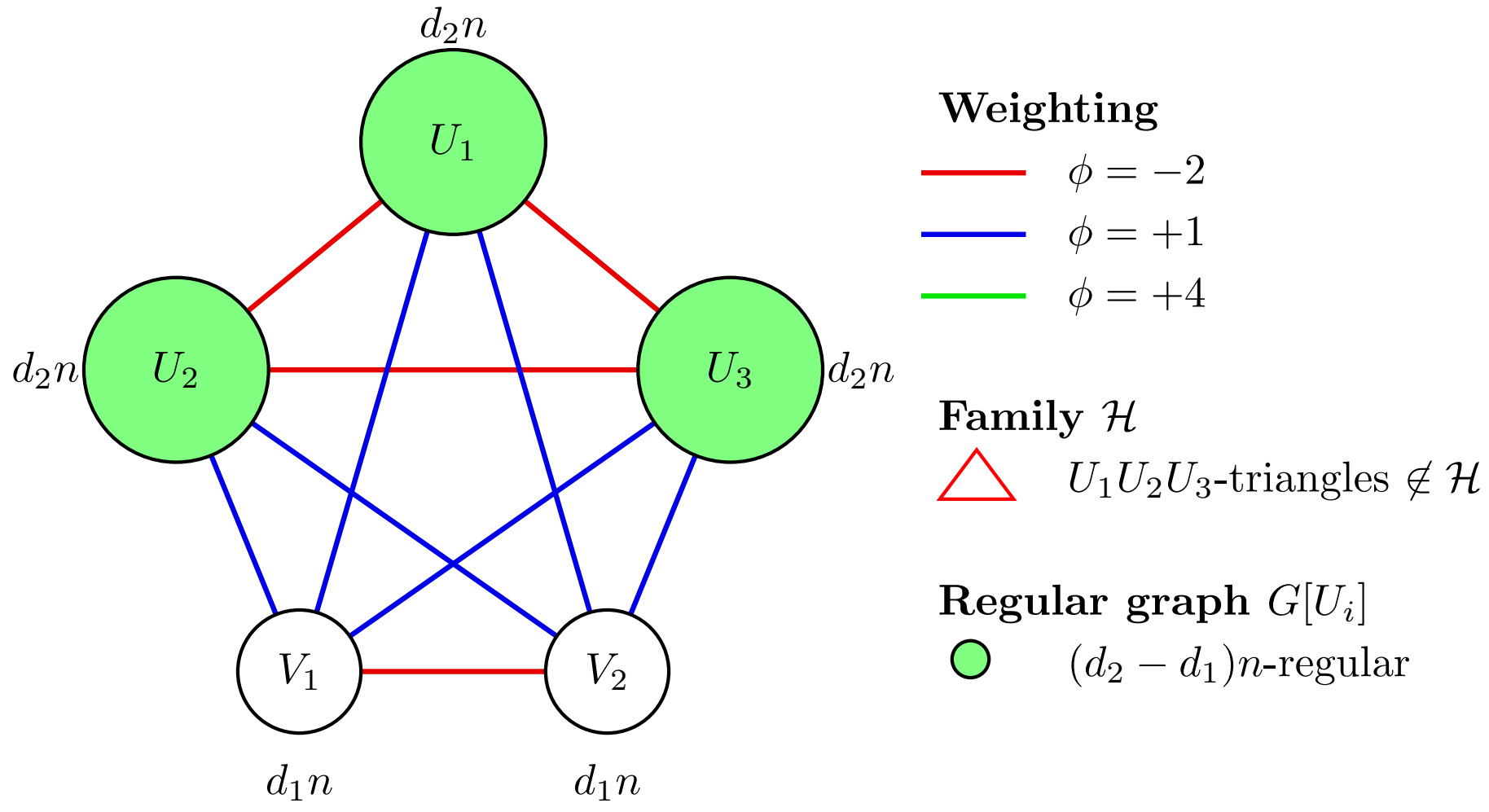

    \centering
    \InterpolConstruction    
    \caption{Interpolating the previous constructions}
    \label{fig:InterpolConst}
\end{figure}

It remains to check that $(G,\mc{H})$ and $\phi$ have the desired properties. It is easy to see that $\phi$ is $\mc{H}$ covered, as every triangle in $\mc{H}$ is either strictly positive (inside one part $U_i$), or sums to $0$. 
We now check that $\phi$ is totally-negative. Indeed we have,
    \[|\phi|=(-2)\cdot \brackets*{3\cdot(d_2n)^2+(d_1n)^2} 
    + (2d_1n)\cdot(3d_2n)+ 4\cdot 3\cdot \frac{(d_2-d_1)^+\cdot n\cdot d_2n}{2}.\]
If $d_1\geq d_2$ we obtain
    \[|\phi|
    =\brackets*{6d_1d_2-6d_2^2-2d_1^2}\cdot n^2
    =\brackets*{-6\cdot\parens*{d_2-\frac{d_1}{2}}-\frac{d_1^2}{2}}\cdot n^2
    <0.
    \]

If $d_2>d_1$  we obtain
    \[|\phi|
    =\brackets*{6d_1d_2-6d_2^2-2d_1^2 + 6d_2^2 - 6d_1d_2}\cdot n^2
    =-2d_1^2\cdot n^2
    <0,
    \]
as $d_1>0$. Therefore $|\phi|<0$, and by~\cref{lem:Farkas}, $(G,\mc{H})$ admits no fractional triangle-decomposition. As we have 
    \[\delta(G)=(1-d_1)n,\quad \Delta^{\rm miss}_2(\mc{H},G)=d_2n,\]
we obtain that $(d_1,d_2)\not\in\mc{D}^*$, and therefore every point above the line $2d_1+3d_2=1$ is also trivially not in $\mc{D}^*$, as desired.
\end{proof}

\section{Concluding Remarks}\label{sec:Conc}
\subsection{Second order discharging}

The sharp fractional Nash-Williams theorem of Delcourt and Postle~\cite{DP26NW} cannot be obtained from first-order proportional discharging alone. Their proof uses a second-order, demand-based rule together with strict sponsorship. Roughly speaking, each negative edge demands charge from its positive neighbors, and each positive edge distributes its available charge in proportion to these demands. The key difficulty in their proof is then to control the total demand at one positive edge $uv$, which is based on properties the common neighborhood of $u,v$, which they call {\em triangle-thinness}.

A direct equivalent program can be derived in our setting by allowing demands only through triangles belonging to $\mc{H}$. The argument that ensures non-negative final charge on positive edges easily survives this restriction. However, the triangle-thinness argument, used to demonstrate that every negative edge receives at least its demand, is strongly impacted by the missing codegree of $\mc{H}$. We conjecture that a second-order, sponsorship-based
extension of the present framework should prove
\[\mc{D}_{\frac14,\frac13}=\set*{(d_1,d_2)\in(0,1]\times[0,1]\colon 4d_1+3d_2\leq1}\subseteq \mc{D}^*.\]

\subsection{Conjectures and potential further progress}

Even if a second-order discharging rule were to yield $\mc{D}_{\frac14,\frac13}\subseteq\mc{D}^*$, this might not yet be the true shape of $\mc{D}^*$. Although the parameter $d_1+d_2$ arise naturally in our argument, there is no compelling reason why the border should be linear between the two points $(0,\frac13)$ and $(\frac14,0)$. In fact, as explained in~\cref{sec:IntuitionLower}, we have reasons to believe that the line $3d_2+2d_1=1$ is asymptotically optimal when $d_1$ is close to $0$, not the current line $3d_2+4d_1=1$ which is the frontier of $\mc{D}_{\frac14,\frac13}$.

\begin{figure}[H]
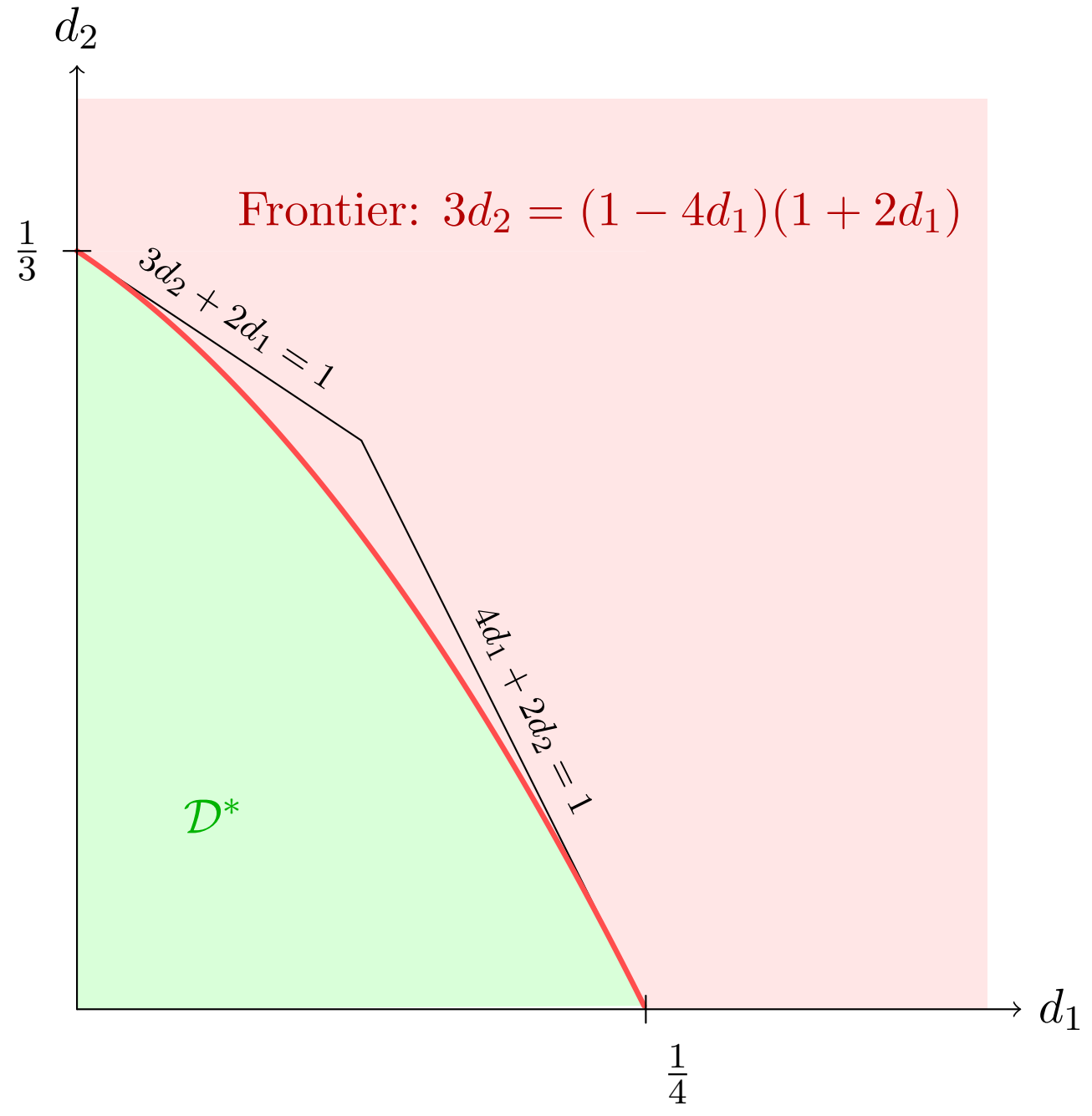

    \centering
    \WildConjecture{18}{0.21}{0.27}
    \caption{One possible conjecture for $\mc{D}^*$}
    \label{fig:wild}
\end{figure}

A natural first idea towards getting a more precise estimate on $\mc{D}^*$ would be to find a line of the form $pd_2+4d_1=1$, for some $p\in[1,3]$ which is asymptotically optimal when $d_2$ is close to $0$. Then one could try to prove that any point ``right'' of that line is not in $\mc{D}^*$. A natural conjecture would be either that the resulting domain is best possible, or that the frontier is defined by some curve with tangents $3d_2+2d_1=1$ when $d_1\to0$ and $pd_2+4d_1=1$ when $d_2\to0$, for example by the curve
\[d_2=\frac{1-4d_1}{3}\cdot\parens*{1+2d_1+\frac{24(2-p)}{p}\cdot d_1^2},\]
which is the unique polynomial of degree at most $3$ with the suggested tangent lines. The case $p=2$ would yield a particularly pleasant quadratic frontier, see~\cref{fig:wild}. The case $p=1$ is also of interest, as it interpolates between the standard Nash-Williams' construction, and the point $d_1=d_2=\frac15$ which, on the line $3d_2+2d_1=1$, corresponds exactly to the balanced complete-$5$-partite construction in~\cref{fig:InterpolConst} with empty parts $U_i$.

\paragraph{AI Declaration.} ChatGPT~5.6 has been used to improve some of the Tikz drawings in this articles. All ideas, proofs, and general writing were developed entirely by the authors in the traditional style.

\bibliographystyle{plain}
\bibliography{bibliography}

\end{document}

%% file: tikz.tex
\usepackage{tikz}
\usepackage{calc}
\usetikzlibrary{calc,backgrounds,decorations.pathreplacing}
\usepackage{ifthen}

\tikzstyle{vertex}=[circle, draw, fill=black, inner sep=0pt, minimum size=4pt]
\newcommand{\vertex}{\node[vertex]}

\newcommand{\DischargingQuasiNeg}[1]
{\begin{tikzpicture}[x=#1 cm, y=#1 cm]
    \vertex (a) at (-30:1) [label=0:$v$]{};
    \vertex (b) at (90:1) [label=90:$u$]{};
    \vertex (c) at (210:1) [label=180:$w$]{};

    \path 
        (a) edge[color=blue] node[below]{$b$} (c)
        (a) edge[color=red] node[right]{$r_2$} (b)
        (c) edge[color=red] node[left]{$r_1$} (b);       

\end{tikzpicture}}

\newcommand{\DischargingQuasiPos}[1]
{\begin{tikzpicture}[x=#1 cm, y=#1 cm]
    \vertex (a) at (-30:1) [label=0:$v$]{};
    \vertex (b) at (90:1) [label=90:$u$]{};
    \vertex (c) at (210:1) [label=180:$w$]{};

    \path 
        (a) edge[color=blue] node[below]{$b$} (c)
        (a) edge[color=red] node[right]{$r$} (b)
        (c) edge[color=blue] node[left]{} (b);       

\end{tikzpicture}}

\newcommand{\RedEdges}[1]
{\begin{tikzpicture}[x=#1 cm, y=#1 cm]

    \vertex (x) at (0,0) [label=0:$u$]{};
    \vertex (y) at (-1,0) [label=200:$v$]{};
    \path (x) edge[color=red, line width=1.5] node[below]{$r$} (y);
    
    \vertex(lf1) at (-2,1)[]{};
    \path (y) edge[color=blue]  (lf1)
          (x) edge[color=red]  (lf1);
    \begin{scope}[on background layer]
        \draw[fill=gray!50, opacity=0.75] (-2,1) ellipse (0.35 and 0.1) node[above=6pt]{$R_u$};
    \end{scope}

    \vertex(rf1) at (1,1)[]{};
    \path (y) edge[color=red]  (rf1)
          (x) edge[color=blue]  (rf1);
    \begin{scope}[on background layer]
        \draw[fill=gray!50, opacity=0.75] (1,1) ellipse (0.35 and 0.1) node[above=6pt]{$R_v$};
    \end{scope}

    \vertex(sk1) at (-0.5,1)[]{};
    \path (y) edge[color=blue]  (sk1)
          (x) edge[color=blue]  (sk1);
    \begin{scope}[on background layer]
        \draw[fill=gray!50, opacity=0.75] (sk1) ellipse (0.35 and 0.1) node[above=6pt]{$B_{uv}$};
    \end{scope}

    \begin{scope}[shift={(0.75,0)}]
        \vertex(d1) at (-0.8,-1)[]{};
        \vertex(d2) at (-0.6,-1)[]{};
        \vertex(d3) at (-0.4,-1)[]{};
        \vertex(d4) at (-0.2,-1)[]{};    
        
        \path (y) edge[color=red]  (d1)
              (x) edge[color=red]  (d1);
        \path (y) edge[color=blue] (d2)
              (x) edge[color=red]  (d2);          
        \path (y) edge[color=red]  (d3)
              (x) edge[color=blue] (d3);          
        \path (y) edge[color=blue] (d4)
              (x) edge[color=blue] (d4);          
        \begin{scope}[on background layer]
            \draw[fill=gray!50, opacity=0.75] (-0.5,-1) ellipse (0.7 and 0.1) node[below=6pt]{$N_{\mc{H}}$};
        \end{scope}

    \end{scope}

    \begin{scope}[shift={(-1,0)}]
        \vertex(e1) at (-0.8,-1)[]{};
        \vertex(e2) at (-0.6,-1)[]{};
        \vertex(e3) at (-0.4,-1)[]{};
        \vertex(e4) at (-0.2,-1)[]{};
        \vertex(e5) at (0,-1)[]{};    
        
        \path (x) edge[color=red]  (e5);
        \path (x) edge[color=blue] (e4);          
        \path (y) edge[color=red]  (e2);          
        \path (y) edge[color=blue] (e1);          
        \begin{scope}[on background layer]
            \draw[fill=gray!50, opacity=0.75] (-0.4,-1) ellipse (0.7 and 0.1) node[below=6pt]{$N_G$};
        \end{scope}

    \end{scope}    

    \draw [decorate,decoration={brace,amplitude=5pt}] (1.5,1.1) -- (1.5,0.05) node[midway,right,xshift=5pt]{$uvw\in\mc{H}$};
    \draw [decorate,decoration={brace,amplitude=5pt}] (1.5,-0.05) -- (1.5,-1.1) node[midway,right,xshift=5pt]{$uvw\not\in\mc{H}$};

\end{tikzpicture}}

\tikzstyle{disc}=[circle, draw, fill=none, inner sep=0, minimum size=60pt,transform shape] 

\newcommand{\LeeLowerConstruction}[2]
{
\begin{tikzpicture}[scale=#1]

  \coordinate (V0) at (0,4);
  \coordinate (V1) at (0,0);
  \coordinate (V2) at (4,0);
  \coordinate (V3) at (4,4);

  \begin{scope}[every label/.style={transform shape}]
      \node[disc, label={[label distance=-0.6cm]90:$V_0$}] (A) at (V0) {};
      \node[disc, label={[label distance=-0.6cm]90:$V_1$}] (B) at (V1) {};
      \node[disc, label={[label distance=-0.6cm]90:$V_2$}] (C) at (V2) {};
      \node[disc, label={[label distance=-0.6cm]90:$V_3$}] (D) at (V3) {};
  \end{scope}

  \node[transform shape] at (2,5) {$E_{#2}$};

  \ifthenelse{#2=0}{
  \begin{scope}[shift={(0,4)}]
    \fill[draw=blue!50, fill=blue!10] (10:1.2) -- (-10:1.2) -- (0:2.8) -- cycle;
    \fill[draw=blue!50, fill=blue!10] (-80:1.2) -- (-100:1.2) -- (-90:2.8) -- cycle;
    \fill[draw=blue!50, fill=blue!10] (-35:1.2) -- (-55:1.2) -- (-45:4.5) -- cycle;
  \end{scope}}{}

  \ifthenelse{#2=1}{
  \begin{scope}[shift={(0,0)}]
    \fill[draw=blue!50, fill=blue!10] ([shift={(0,4)}] -35:1.2) -- ([shift={(4,4)}]-145:1.2) -- (35:1.2) -- cycle;
    \begin{scope}[on background layer]
        \fill[draw=blue!50, fill=blue!10] ([shift={(0,4)}] -45:1.2) -- ([shift={(4,4)}]-135:1.2) -- ([shift={(4,0)}] 135:1.2) -- cycle;  
        \fill[draw=blue!50, fill=blue!10] ([shift={(0,4)}] -50:1.2) -- ([shift={(0,0)}]46:1.10) -- ([shift={(4,0)}] 140:1.2) -- cycle;          
    \end{scope}
  \end{scope}}{}  

  \ifthenelse{#2=2}{
  \newcommand{\bluetriangle}{\fill[draw=blue!50, fill=blue!10] (90:0.6) -- (-30:0.6) -- (-150:0.6) -- cycle;}
  \begin{scope}[shift={(0,-0.2)}]
    \bluetriangle
  \end{scope}
  \begin{scope}[shift={(4,3.8)}]
    \bluetriangle
  \end{scope}
  \begin{scope}[shift={(4,-0.2)}]
    \bluetriangle
  \end{scope}
  }{}

  \ifthenelse{#2=3}{
  \begin{scope}
    \fill[draw=blue!50, fill=blue!10] (40:1.2) -- (60:1.2) -- ([shift={(4,4)}]-150:1.2) -- cycle;
    \fill[draw=blue!50, fill=blue!10] (30:1.2) -- ([shift={(4,4)}] -140:1.2) -- ([shift={(4,4)}]-120:1.2) -- cycle;

    \fill[draw=blue!50, fill=blue!10] ([shift={(0,0.2)}] -10:1.2) -- ([shift={(0,0.2)}] 10:1.2) -- ([shift={(4,0.2)}]180:1.2) -- cycle;
    \fill[draw=blue!50, fill=blue!10] ([shift={(4,-0.2)}] 170:1.2) -- ([shift={(4,-0.2)}] 190:1.2) -- ([shift={(0,-0.2)}]0:1.2) -- cycle;

    \fill[draw=blue!50, fill=blue!10] ([shift={(3.8,4)}]-100:1.2) -- ([shift={(3.8,4)}]-80:1.2) -- ([shift={(3.8,0)}]90:1.2) -- cycle;
    \fill[draw=blue!50, fill=blue!10] ([shift={(4.2,0)}]100:1.2) -- ([shift={(4.2,0)}]80:1.2) -- ([shift={(4.2,4)}]-90:1.2) -- cycle;
    
  \end{scope}}{}  

\end{tikzpicture}}

\newcommand{\LeeBis}{
\begin{tikzpicture}[
    scale=1,
    every node/.style={font=\small},
    part/.style={
        circle,
        draw=blue,
        thick,
        fill=blue!25!white,
        minimum size=13mm,
        inner sep=1pt
    },
    negative/.style={draw=red!75!black, very thick},
    positive/.style={draw=blue!70!black, very thick},
    forbidden/.style={draw=gray!70!black, dashed, thick}
]

\coordinate (cV1) at (0,2.3);
\coordinate (cV2) at (-1.8,-0.7);
\coordinate (cV3) at (1.8,-0.7);

\draw[negative] (cV1) -- (cV2);
\draw[negative] (cV1) -- (cV3);
\draw[negative] (cV2) -- (cV3);

\node[part] (V1) at (cV1) {$U_1$};
\node[part] (V2) at (cV2) {$U_2$};
\node[part] (V3) at (cV3) {$U_3$};

\node at (0,0) {$K_n$};

\node at (0,3.15) {$n/3$};
\node[anchor=east] at (-2.55,-0.7) {$n/3$};
\node[anchor=west] at (2.55,-0.7) {$n/3$};

\node[anchor=west,font=\small\bfseries] at (3.5,2.25)
    {Weighting};

\draw[negative] (3.5,1.75) -- (4.35,1.75);
\node[anchor=west] at (4.55,1.75)
    {$\phi=-1$};

\draw[positive] (3.5,1.25) -- (4.35,1.25);
\node[anchor=west] at (4.55,1.25)
    {$\phi=+2$};

\node[anchor=west,font=\small\bfseries] at (3.5,0.25)
    {Family $\mc{H}$};
\draw[red, thick]
    (3.95,-0.4) -- (4.25,0) -- (4.55,-0.4) -- cycle;
\node[anchor=west, align=left] at (4.85,-0.25)
    {$V_1V_2V_3$-triangles~$\not\in \mc{H}$};

\end{tikzpicture}
}

\newcommand{\DomainConjecture}[3]{%
\begin{tikzpicture}[
    scale=#1,
    line cap=round,
    line join=round
]
    \def\xmax{0.40}
    \def\ymax{0.40}

    \def\xcrit{0.25}          
    \def\ycrit{0.3333333333}  
    \def\yqright{0.1666666667}

    \def\xLee{0.06042}
    \def\yLee{0.12084}
    \def\xZheng{0.071083}
    \def\yZheng{0.142166}

    \fill[red!20]
        (0,\ycrit) rectangle (\xmax,\ymax);

    \fill[red!20]
        (\xcrit,0) rectangle (\xmax,\ycrit);

    \fill[red!40]
        (0,\ycrit)
        -- (\xcrit,\ycrit)
        -- (\xcrit,\yqright)
        -- cycle;

    \fill[white]
        (0,0)
        -- (0,\ycrit)
        -- (\xcrit,\yqright)
        -- (\xcrit,0)
        -- cycle;

    \fill[green!20]
        (0,0)
        -- (0,\ycrit)
        -- (\xcrit,0)
        -- cycle;

    \fill[blue!20]
        (0,0)
        -- (0,#3)
        -- (#2,0)
        -- cycle;

    \fill[blue!50]
        (0,0)
        -- (0,\yZheng)
        -- (\xZheng,0)
        -- cycle;             

    \fill[blue!80]
        (0,0)
        -- (0,\yLee)
        -- (\xLee,0)
        -- cycle;          

    \draw[thick, red]
        (0,\ycrit)
        -- (\xcrit,\yqright);

    \draw[thick,blue!90!black]
        (0,#3) -- (#2,0);

    \draw (0,0) -- (0.32,0);
    \draw[->] (0.33,0) -- (\xmax+0.015,0)
        node[right] {$d_1$};
    
    \draw[thick]
        (0.32,-0.008) -- (0.325,0.008)
        (0.325,-0.008) -- (0.33,0.008);

    \draw (0,0) -- (0,0.36);
    \draw[->] (0,0.37) -- (0,\ymax+0.015)
        node[above] {$d_2$};
    
    \draw[thick]
        (-0.008,0.36) -- (0.008,0.365)
        (-0.008,0.365) -- (0.008,0.37);        

    \draw (\xmax,0.006) -- (\xmax,-0.006) node[anchor=north] {$1$};
    \draw (0.006,\xmax) -- (-0.006,\xmax) node[anchor=east] {$1$};

    \draw (#2,0.006) -- (#2,-0.006) node[anchor=north] {$#2$};
    \draw (0.006,#3) -- (-0.006,#3) node[anchor=east] {$#3$};
        
    \draw (\xcrit,0.006) -- (\xcrit,-0.006) node[anchor=north] {$\frac14$};
    \draw (0.006,\ycrit) -- (-0.006,\ycrit) node[anchor=east] {$\frac13$};

    \draw (\xLee,0.006) -- (\xLee,-0.006) node[anchor=north,xshift=-8pt] {$\approx0.06$};
    \draw (0.006,\yLee) -- (-0.006,\yLee) node[anchor=east] {$\approx0.12$};

    \draw (\xZheng,0.006) -- (\xZheng,-0.006) node[anchor=north,xshift=8pt,yshift=-10pt] {$\approx0.07$};
    \draw (0.006,\yZheng) -- (-0.006,\yZheng) node[anchor=east] {$\approx0.14$};
    
    \draw (0.006,\yqright) -- (-0.006,\yqright) node[anchor=east] {$\frac16$};

    \begin{scope}[shift={(0.46,0.335)}]

        \node[anchor=west,font=\scriptsize] at (-0.01,0.04) {\underline{$\mathbf{(d_1,d_2)\in\mc{D}^*}$}};
        \begin{scope}{shift={(0,0)}}
        \fill[blue!80]
            (0,0) rectangle (0.025,0.02);
        \node[anchor=west,font=\scriptsize]
            at (0.03,0.01)
            {Implied by Lee~\cite{L25}};
        \end{scope}

        \begin{scope}[shift={(0,-0.03)}]
        \fill[blue!50]
            (0,0) rectangle (0.025,0.02);
        \node[anchor=west,font=\scriptsize]
            at (0.03,0.01)
            {Implied by Zheng~\cite{Z24}};
        \end{scope}

        \begin{scope}[shift={(0,-0.06)}]
        \fill[blue!20]
            (0,0) rectangle (0.025,0.02);
        \node[anchor=west,font=\scriptsize]
            at (0.03,0.01)
            {\Cref{thm:MainSimplicial}:
             $\mathcal D_{#2,#3}\subseteq\mathcal D^*$};
        \end{scope}

        \node[anchor=west,font=\scriptsize] at (-0.01,-0.10) {\underline{\textbf{Open question}~}};        
        \begin{scope}[shift={(0,-0.14)}]
        \fill[green!20]
            (0,0) rectangle (0.025,0.02);
        \node[anchor=west,font=\scriptsize]
            at (0.03,0.01)
            {Conjectured:
            $\mathcal D_{\frac14,\frac13}
            \subseteq\mathcal D^*$};
        \end{scope}          

        \begin{scope}[shift={(0,-0.17)}]
        \draw (0,0) rectangle (0.025,0.02);        
        \node[anchor=west,font=\scriptsize]
            at (0.03,0.01)
            {?};
        \end{scope}

        \node[anchor=west,font=\scriptsize] at (-0.01,-0.21) {\underline{$\mathbf{(d_1,d_2)\notin\mc{D}^*}$}}; 
        \begin{scope}[shift={(0,-0.25)}]
        \fill[red!40]
            (0,0) rectangle (0.025,0.02);
        \node[anchor=west,font=\scriptsize]
            at (0.03,0.01)
            {Excluded by~\cref{lowerdomain}};
        \end{scope}     
        
        \begin{scope}[shift={(0,-0.28)}]
        \fill[red!20]
            (0,0) rectangle (0.025,0.02);
        \node[anchor=west,font=\scriptsize]
            at (0.03,0.01)
            {Excluded outside $[0,\frac14]\times[0,\frac13]$};
        \end{scope}

    \end{scope}

\end{tikzpicture}%
}

\newcommand{\FirstConstruction}{
\begin{tikzpicture}[
    scale=1,
    every node/.style={font=\small},
    smallpart/.style={
        circle,
        draw=black,
        thick,
        fill=white,
        minimum size=12mm,
        inner sep=1pt
    },
    largepart/.style={
        circle,
        draw=black,
        thick,
        fill=white,
        minimum size=14mm,
        inner sep=1pt
    },    
    neg/.style={draw=red!90!black, very thick},
    posit/.style={draw=blue!90!black, very thick}
]

\draw[neg] (0,2.5) -- (-2.25,0.65);
\draw[neg] (0,2.5) -- (2.25,0.65);
\draw[neg] (-2.25,0.65) -- (2.25,0.65);
\draw[neg] (-1.25,-1.8) -- (1.25,-1.8);

\draw[posit] (0,2.5) -- (-1.25,-1.8);
\draw[posit] (0,2.5) -- (1.25,-1.8);

\draw[posit] (-2.25,0.65) -- (-1.25,-1.8);
\draw[posit] (-2.25,0.65) -- (1.25,-1.8);

\draw[posit] (2.25,0.65) -- (-1.25,-1.8);
\draw[posit] (2.25,0.65) -- (1.25,-1.8);

\node[smallpart] (U1) at (0,2.5) {$U_1$};
\node[smallpart] (U2) at (-2.25,0.65) {$U_2$};
\node[smallpart] (U3) at (2.25,0.65) {$U_3$};
\node[largepart] (V1) at (-1.25,-1.8) {$V_1$};
\node[largepart] (V2) at (1.25,-1.8) {$V_2$};

\node at (0,3.35) {$n/6$};
\node[anchor=east] at (-2.9,0.65) {$n/6$};
\node[anchor=west] at (2.9,0.65) {$n/6$};

\node at (-1.25,-2.75) {$n/4$};
\node at (1.25,-2.75) {$n/4$};

\node[anchor=west,font=\small\bfseries] at (3.8,2.25)
    {Weighting};

\draw[neg] (3.8,1.75) -- (4.65,1.75);
\node[anchor=west] at (4.85,1.75) {$\phi=-2$};

\draw[posit] (3.8,1.25) -- (4.65,1.25);
\node[anchor=west] at (4.85,1.25) {$\phi=+1$};

\node[anchor=west,font=\small\bfseries] at (3.8,0.25)
    {Family $\mc{H}$};
\draw[red, thick]
    (3.95,-0.4) -- (4.25,0) -- (4.55,-0.4) -- cycle;
\node[anchor=west, align=left] at (4.85,-0.25)
    {$U_1U_2U_3$-triangles~$\not\in \mc{H}$};

\end{tikzpicture}}

\newcommand{\InterpolConstruction}{
\begin{tikzpicture}[
    scale=1,
    every node/.style={font=\small},
    smallpart/.style={
        circle,
        draw=black,
        thick,
        fill=green!50!white,
        minimum size=15mm,
        inner sep=1pt
    },
    largepart/.style={
        circle,
        draw=black,
        thick,
        fill=white,
        minimum size=10mm,
        inner sep=1pt
    },    
    legendpart/.style={
        circle,
        draw=black,
        thick,
        fill=green!50!white,
        minimum size=3mm,
        inner sep=1pt
    },    
    neg/.style={draw=red!90!black, very thick},
    posit/.style={draw=blue!90!black, very thick},
    positbis/.style={draw=green!90!black, very thick}
]

\draw[neg] (0,2.5) -- (-2.25,0.65);
\draw[neg] (0,2.5) -- (2.25,0.65);
\draw[neg] (-2.25,0.65) -- (2.25,0.65);
\draw[neg] (-1.25,-1.8) -- (1.25,-1.8);

\draw[posit] (0,2.5) -- (-1.25,-1.8);
\draw[posit] (0,2.5) -- (1.25,-1.8);

\draw[posit] (-2.25,0.65) -- (-1.25,-1.8);
\draw[posit] (-2.25,0.65) -- (1.25,-1.8);

\draw[posit] (2.25,0.65) -- (-1.25,-1.8);
\draw[posit] (2.25,0.65) -- (1.25,-1.8);

\node[smallpart] (U1) at (0,2.5) {$U_1$};
\node[smallpart] (U2) at (-2.25,0.65) {$U_2$};
\node[smallpart] (U3) at (2.25,0.65) {$U_3$};
\node[largepart] (V1) at (-1.25,-1.8) {$V_1$};
\node[largepart] (V2) at (1.25,-1.8) {$V_2$};

\node at (0,3.45) {$d_2n$};
\node[anchor=east] at (-2.9,0.65) {$d_2n$};
\node[anchor=west] at (2.9,0.65) {$d_2n$};

\node at (-1.25,-2.70) {$d_1n$};
\node at (1.25,-2.70) {$d_1n$};

\node[anchor=west,font=\small\bfseries] at (3.8,2.75)
    {Weighting};

\draw[neg] (3.8,2.25) -- (4.65,2.25);
\node[anchor=west] at (4.85,2.25) {$\phi=-2$};

\draw[posit] (3.8,1.75) -- (4.65,1.75);
\node[anchor=west] at (4.85,1.75) {$\phi=+1$};

\draw[positbis] (3.8,1.25) -- (4.65,1.25);
\node[anchor=west] at (4.85,1.25) {$\phi=+4$};

\node[anchor=west,font=\small\bfseries] at (3.8,0.25)
    {Family $\mc{H}$};
\draw[red, thick]
    (3.95,-0.4) -- (4.25,0) -- (4.55,-0.4) -- cycle;
\node[anchor=west, align=left] at (4.85,-0.25)
    {$U_1U_2U_3$-triangles~$\not\in \mc{H}$};

\node[anchor=west,font=\small\bfseries] at (3.8,-1.25)
    {Regular graph $G[U_i]$}; 
\node[legendpart] (legpart) at (4.2,-1.70) {};
\node[anchor=west, align=left] at (4.85,-1.75)
    {$(d_2-d_1)n$-regular};
    
\end{tikzpicture}}

\newcommand{\WildConjecture}[3]{%
\begin{tikzpicture}[
    scale=#1,
    line cap=round,
    line join=round
]
    \def\xmax{0.40}
    \def\ymax{0.40}

    \def\xcrit{0.25}          
    \def\ycrit{0.3333333333}  
    
    \def\xmid{0.125}
    \def\ymid{0.25}

    \fill[red!10]
        (0,\ycrit) rectangle (\xmax,\ymax);

    \fill[red!10]
        (\xcrit,0) rectangle (\xmax,\ycrit);

    \fill[green!15]
        (0,0)
        -- (0,1/3)
        -- plot[domain=0:0.25, samples=100]
            (\x,{(1-4*\x)*(1+2*\x)/3})
        -- cycle;

    \fill[red!10]
        (0,1/3)
        -- (0.25,1/3)
        -- (0.25,0)
        -- plot[domain=0.25:0, samples=100]
            (\x,{(1-4*\x)*(1+2*\x)/3})
        -- cycle;

    \draw[black]
        (0,\ycrit)
        -- (\xmid,\ymid)
        node[midway, above, sloped, font=\scriptsize]
    {$3d_2+2d_1=1$};

    \draw[black]
        (\xcrit,0)
        -- (\xmid,\ymid)
        node[midway, above, sloped, font=\scriptsize]
    {$4d_1+2d_2=1$};        

    \draw[red!70, very thick, domain=0:0.25, samples=100, smooth]
        plot (\x, {(1-4*\x)*(1+2*\x)/3});

    \draw[->]
        (0,0) -- (\xmax+0.015,0)
        node[right] {$d_1$};

    \draw[->]
        (0,0) -- (0,\ymax+0.015)
        node[above] {$d_2$};

    \draw
        (\xcrit,0.006) -- (\xcrit,-0.006)
        node[anchor=north west,yshift=-1pt] {$\frac14$};

    \draw
        (0.006,\ycrit) -- (-0.006,\ycrit)
        node[anchor=east,xshift=-1pt] {$\frac13$};

    \node[
        green!70!black
        ] at (0.060,0.085)
        {$\mathcal{D}^*$};

    \node[
        red!70!black
    ] at (0.23,0.35)
        {Frontier: $3d_2=(1-4d_1)(1+2d_1)$};

\end{tikzpicture}%
}